\RequirePackage{fix-cm}
\documentclass[smallextended]{svjour3}
\smartqed  % flush right qed marks, e.g. at end of proof
\usepackage[hyphens]{url}
\usepackage{hyperref}
\usepackage[hyphenbreaks]{breakurl}
\usepackage{graphicx}
\usepackage{bm}
\usepackage{tikz}
\usepackage{pgfplots}
\pgfplotsset{compat=1.18}
\usetikzlibrary{plotmarks}
\usepackage[english]{babel}
\usepackage{times}
\usepackage{orcidlink}
\usepackage{listings}
\usepackage{booktabs}
\usepackage{multirow}
\usepackage{amssymb}
\usepackage{amsmath}
\usepackage{setspace}
\usepackage{blindtext}
\usepackage{xspace}
\usepackage{algpseudocode}
\usepackage{algorithm}
\usepackage{csquotes}
\usepackage{comment}
\usepackage{enumitem}

\begin{document}

\setlength{\abovecaptionskip}{-12pt}
\setlength{\belowcaptionskip}{12pt}
\setlength{\floatsep}{0pt}
\setlength{\textfloatsep}{10pt plus 0pt minus 2.0pt}
\setlength{\intextsep}{10pt plus 0pt minus 2.0pt}

\title{Time-limited Balanced Truncation for Data Assimilation Problems
\thanks{
\begin{acknowledgements}
We thank Patrick Kürschner for providing his code for time-limited balanced truncation and helpful insights on this topic. We also thank Elizabeth Qian for the fruitful discussions about her work and during her visit to the University of Potsdam. We thank the reviewers and Thomas Mach for critically reading the manuscript and suggesting substantial improvements. Our research has been partially funded by the Deutsche Forschungsgemeinschaft (DFG) --- Project-ID 318763901 -- SFB1294.
\end{acknowledgements}}}

\author{Josie König \orcidlink{0000-0003-4999-6649} \and Melina A.\@ Freitag \orcidlink{0000-0002-4539-2162}}

\institute{Josie König \at
            Universität Potsdam, Institut für Mathematik\\
            Karl-Liebknecht-Str. 24--25, D--14476 Potsdam \\
            Tel.: +49 331 977--230182\\
            \email{josie.koenig@uni-potsdam.de}
            \and Melina A.\@ Freitag \at
            Universität Potsdam, Institut für Mathematik\\ \email{melina.freitag@uni-potsdam.de}}
            
\date{Version as of 24.08.2023}
\maketitle
\begin{abstract}
Balanced truncation is a well-established model order reduction method which has been applied to a variety of problems. Recently, a connection between linear Gaussian Bayesian inference problems and the system-theoretic concept of balanced truncation has been drawn \cite{Qian2021Balancing}. Although this connection is new, the appli\-cation of balanced truncation to data assimilation is not a novel idea: it has already been used in four-dimensional variational data assimilation (4D-Var). This paper dis\-cusses the application of balanced truncation to linear Gaussian Bayesian inference, and, in particular, the 4D-Var method, thereby strengthening the link between systems theory and data assimilation further. Similarities between both types of data \linebreak assimilation problems enable a generalisation of the state-of-the-art approach to the use of arbitrary prior covariances as reachability Gramians. Furthermore, we propose an enhanced approach using time-limited balanced truncation that allows to balance Bayesian inference for unstable systems and in addition improves the numerical results for short observation periods.
\end{abstract}

\section{Introduction}
One of the greatest challenges of 21st century mathematics is integrating large a\-mounts of data into computational models to obtain novel insights into the state and dynamics of a system. This is the typical framework for data assimilation. Data assimilation combines available, mostly high-dimensional, models with observations at certain time steps, often through time-dependent PDEs. These models occur in various domains of application, such as meteorology \cite{Dee1991Meteor,Ghil1989Meteor,Ghil1991Meteor}, geosciences \linebreak\mbox{\cite{Carr2018Geosc,Fletcher2017GeoSc,Park2017GeoSc}}, and medicine \cite{Engbert2020Med,Kost2011Med}. They are often unstable.

The typical approach to data assimilation problems is Bayesian inference; see, for instance, \cite{Dashti2017BIPSurvey,Law2015DA,Taran2005BIPbook}. We consider the parameters $\pmb{x} \in \mathbb{R}^d$ to be the realisations of a random variable. The distribution of said random variable is what we want to infer. Some \textit{prior} information is given. The time-dependent observations $\pmb{m}$ are related to the parameters via a given forward model operator $\pmb{G}$ with additive noise $\pmb{\epsilon}$:
\begin{align}\label{DynSys6}
    \pmb{m} = \pmb{G}(\pmb{x})+ \pmb{\epsilon}.
\end{align}
From \eqref{DynSys6}, the likelihood $\pmb{m}|\pmb{x}$ is derived. Combining this with the measurements, the prior is updated to the posterior distribution of interest $\pmb{x}|\pmb{m}$. Usually, the initial state of a dynamical system, or rather its distribution, is to be inferred.

Often, the state is high-dimensional and therefore forward model simulations are expensive. This problem is well known from systems and control theory, where reduced-order modelling has been established; see, e.g., \cite{Antoulas2005Book,Benner2015Survey,Benner2017MORBook}. Both non-intrusive data-based methods \cite{Benner2020NonInt,Pehers2016NonInt} and intrusive model-based methods \cite{Antoulas2020interpolatory,Gugercin2004Survey} have been in the spotlight of past and current research. Model-based approaches are promising in the data assimilation context. One such method is balanced truncation (BT) reducing a stable, linear time-invariant system by projecting it onto a subspace of simul\-ta\-neous\-ly easily reachable and observable states. Variants for time-varying systems \mbox{\cite{Lang2016TimeVar,Sandberg2004TimeVar}}, inhomogeneous systems \mbox{\cite{Beattie2016inhom,Heink2011inhom}}, bilinear systems \cite{Benner2017BilinBT,Duff2019BilinBT} and even nonlinear systems \cite{kramer2022nonlinear,kramer2022nonlinear2} have been developed. The application of BT to data assimilation has appeared first in the context of 4D-Var \cite{BoessDiss,GreenDiss,Lawless20084D-Varproposal,Lawless2008Strong4DVarBT}.

Recently, Qian et al.\@ \cite{Qian2021Balancing} have reinvigorated interest in the connection between system-theoretic model order reduction and Bayesian inference. Balanced truncation for linear Gaussian Bayesian inference has been considered. In this paper, we relate the results of Qian et al. to the application of BT to 4D-Var. We generalise the existing results on BT for linear Gaussian Bayesian inference to arbitrary Gaussian priors and unstable systems by proposing time-limited balanced truncation (TLBT) for Bayesian inference. The proposed method improves on the state-of-the-art approach, especially for small observation intervals appearing in many data assimilation applications.

The remainder of the paper is structured as follows: At the beginning, the frame\-works for two data assimilation problems and the basics of BT are presented; Sect.\@ \ref{Basics} Sect.\@ \ref{BT for DA} introduces how time-limited balanced truncation is applied in these situations. Sect.\@ \ref{TLBT-Exps} presents numerical results before the paper is concluded in Sect.\@ \ref{Summary} with the discussion of directions for future work.

\section{Setting and Background}\label{Basics}
We consider the linear dynamical system
\begin{subequations}\label{DynSys1}
  \begin{align}\label{BLIP_forward}
  \frac{\mathrm{d}\pmb{x}}{\mathrm{d}t}  = \pmb{Ax},
  \end{align}
 where $\pmb{x}(t)\in\mathbb{R}^d$ is the state of the system at time $t \in \mathbb{R}$, and $\pmb{A} \in \mathbb{R}^{d\times d}$ is a linear model operator describing the state evolution. The initial state $\pmb{x}(0) = \pmb{x}_0$ is unknown. Noise-polluted outputs $\pmb{m}_k\in \mathbb{R}^{d_\mathrm{out}}$ of the system are only available at discrete time points $0 < t_1 < \ldots < t_n$ via the measurement model
  \begin{align}\label{GaussMeas}
  \pmb{m}_k = \pmb{C}\pmb{x}(t_k) + \pmb{\epsilon}_k = \pmb{C}e^{\pmb{A}t_k}\pmb{x}_0 + \pmb{\epsilon}_k,
  \end{align}
  with  $\pmb{C} \in \mathbb{R}^{d_\mathrm{out}\times d}$ the state-to-output operator, typically with $d_\mathrm{out} \ll d$. The noise $\pmb{\epsilon}_k \sim \mathcal{N}(\pmb{0},\pmb{\Gamma}_\epsilon)$ is independently and identically distributed additive Gaussian with known positive definite covariance $\pmb{\Gamma}_\epsilon \in \mathbb{R}^{d_\mathrm{out}\times d_\mathrm{out}}$. The observation interval of the system is $[0,t_n]$. We always start at time 0, with the unknown initial condition, and set the end time to $t_e = t_n$.
  
The data assimilation task is to infer the unknown initial state from the noisy measurements. To this end, the Bayesian statistical approach is applied and the initial condition is assumed to be a-priori Gaussian distributed:
\begin{align}\label{GaussPrior}
    \pmb{x}_0 \sim \mathcal{N}(\pmb{0},\pmb{\Gamma}_{\mathrm{pr}}),
\end{align}
with $\pmb{\Gamma}_{\mathrm{pr}} \in \mathbb{R}^{d\times d}$ a given positive definite prior covariance  matrix.
\end{subequations}
    
Substituting the Gaussian prior \eqref{GaussPrior} in the linear Gaussian measurement model \eqref{GaussMeas} yields a Gaussian likelihood for the measurements conditioned on the initial condition and the measurement times $\pmb{t} = [t_1,\ldots,t_n]^\mathrm{T}$:
    \begin{align*}
        \pmb{m}|(\pmb{x}_0,\pmb{t}) \sim \mathcal{N}(\pmb{G}\pmb{x}_0, \pmb{\Gamma}_{\mathrm{obs}}),
    \end{align*}
    where $\pmb{m} \in \mathbb{R}^{nd_\mathrm{out}}$, $\pmb{G} \in \mathbb{R}^{nd_\mathrm{out} \times d}$, and $\pmb{\Gamma}_{\mathrm{obs}} \in \mathbb{R}^{nd_\mathrm{out} \times nd_\mathrm{out}}$ as follows:
    \begin{align*}
        \pmb{m}= \begin{bmatrix}\pmb{m}_1\\\vdots\\\pmb{m}_n
\end{bmatrix},\enspace \pmb{G}= \begin{bmatrix}\pmb{C}e^{\pmb{A}t_1}\\\vdots\\\pmb{C}e^{\pmb{A}t_n}\end{bmatrix}, \enspace \pmb{\Gamma}_{\mathrm{obs}}= \begin{bmatrix}\pmb{\Gamma}_{\pmb{\epsilon}}&&&\\&\pmb{\Gamma}_{\pmb{\epsilon}}&&\\ && \ddots &\\ &&&\pmb{\Gamma}_{\pmb{\epsilon}}
\end{bmatrix}.
    \end{align*}
For a Gaussian prior and a Gaussian likelihood, the posterior is again Gaussian \cite{Carlin2009BayesMeth}:
    \begin{align}\label{PostStats}
        \pmb{x}_0|(\pmb{m},\pmb{t}) &\sim \mathcal{N}(\pmb{\mu}_{\mathrm{pos}}, \pmb{\Gamma}_{\mathrm{pos}}), \enspace \text{where}\notag \\
        \pmb{\mu}_{\mathrm{pos}} &= \pmb{\Gamma}_{\mathrm{pos}}\pmb{G}^{\mathrm{T}}\pmb{\Gamma}_{\mathrm{obs}}^{-1}\pmb{m}\enspace \text{and} \enspace \pmb{\Gamma}_{\mathrm{pos}} = (\pmb{H} + \pmb{\Gamma}_{\mathrm{pr}}^{-1})^{-1}.
    \end{align}
    The matrix $\pmb{H} \in \mathbb{R}^{d \times d}$ is called the Fisher information of $\pmb{m}$ \cite{LehmCase1998PointEst} and defined as
    \begin{align*}
        \pmb{H} = \pmb{G}^{\mathrm{T}}\pmb{\Gamma}_{\mathrm{obs}}^{-1}\pmb{G} = \sum_{k=1}^n \pmb{G}_k^{\mathrm{T}}\pmb{\Gamma}_{\pmb{\epsilon}}^{-1}\pmb{G}_k = \sum_{k=1}^n e^{\pmb{A}^{\mathrm{T}}t_k}\pmb{C}^{\mathrm{T}}\pmb{\Gamma}_{\pmb{\epsilon}}^{-1}\pmb{C}e^{\pmb{A}t_k}.
    \end{align*}

The advantage of this data assimilation setting is that the posterior statistics \eqref{PostStats} can be obtained in a closed form. The challenge is the computation of the posterior mean and covariance for a high-dimensional state space, that is for $d$ large. This implies several multiplications with the high-dimensional forward map $\pmb{G}$ and its transpose, which might only be available through costly evolution of the underly\-ing dynamical system. Benner, Qiu, and Stoll have proposed a low-rank compression of the posterior covariance \cite{Benner2018LRCov}. Other than this, the usually small dimension $d_\mathrm{out}$ of the observations suggests that the dynamics of interest may be modelled accurately when restricted to a low-dimensional subspace. Model order reduction approaches, such as those presented in this paper, aim to reduce the dimensionality of the problem to make computations more efficient while maintaining sufficient accuracy.

\subsection{Incremental Four Dimensional Variational Data Assimilation (4D-Var)}\label{4DVar}
4D-Var is a variational data assimilation technique originating in weather fore\-cast\-ing~\cite{Sasaki70Basic}. In 4D-Var, instead of a continuous model evolution, a discrete time-varying system for $k = 0,1,\ldots,n$ is considered:
\begin{subequations}\label{DynSys3}
\begin{align}
\pmb{x}_{k+1} &= \mathcal{A}_k(\pmb{x}_{k}), \label{DynSys3a}\\
\pmb{m}_k &= \mathcal{C}_k(\pmb{x}_k) + \pmb{\epsilon}_k, \label{DynSys3b}
\end{align}
where $\pmb{x}_k \in \mathbb{R}^d$ is the state of the system at time $t_k \in \{0,\ldots,t_n\}$, and $\mathcal{A}_k:\mathbb{R}^d\rightarrow \mathbb{R}^d$ is the nonlinear model operator that describes the evolution of the state from time $t_k$ to $t_{k+1}$. The measurement model suffers from centered Gaussian observational noise $\pmb{\epsilon}_k\sim \mathcal{N}(\pmb{0},\pmb{R}_k)$ for $k = 0,\ldots,n$.
The goal is again to infer the initial state $\pmb{x}_0$ from the noisy measurements $\{ \pmb{m}_k\}_{k = 0,\ldots,n}$. We assume that an a-priori estimate $\pmb{x}_0^b$ of the initial state is given---the so-called \textit{background state}---which is error-prone:
\begin{align}
    \pmb{x}_0 - \pmb{x}_0^b = \pmb{e}_0 \sim \mathcal{N}(\pmb{0},\pmb{\Gamma}_\mathrm{pr}).
\end{align}
\end{subequations}
Linear Gaussian Bayesian inference and 4D-Var are closely related: discretising~\eqref{DynSys1} with equispaced points in time $t_k = k\cdot h,\, k=0,\ldots,n$ and $t_n = t_e$, we obtain a system of the more general type \eqref{DynSys3} with $\mathcal{A}_k = \pmb{A}_{\mathrm{disc}} = e^{\pmb{A}h}$ and $\mathcal{C}_k = \pmb{C}$ for all~$k$. A closed form solution to the posterior statistics for the initial state can only be computed for systems of type \eqref{DynSys3} if $\mathcal{A}_k$ is linear and time-invariant. In 4D-Var, for arbitrary $\mathcal{A}_k$, the maximum a posteriori estimate for $\pmb{x}_0$, i.e., \mbox{\@ $\pmb{x}_0^* = \arg \min_{\pmb{x}_0} J(\pmb{x}_0)$,} is obtained by minimising the cost functional derived from the posterior distribution; see, for instance, \cite{Melina2020GAMM}:
\begin{subequations}\label{4D-Var_Cost}
\begin{align}
J(\pmb{x}_0) = \frac{1}{2}\left\| \pmb{x}_0 - \pmb{x}_0^b\right\|^2_{\pmb{\Gamma}_\mathrm{pr}^{-1}}&+\frac{1}{2}\sum_{k=0}^n\left\| \pmb{m}_k - \mathcal{C}_k(\pmb{x}_k) \right\|^2_{\pmb{R}_k^{-1}}
\end{align}
subject to the forward model dynamics
\begin{align}
\pmb{x}_{k+1} &= \mathcal{A}_k(\pmb{x}_{k}) \enspace \text{for } k = 0,\ldots, n-1,
\end{align}
\end{subequations}
where $\pmb{x}_k$ is the state at time step $t_k$ and $\left\|\pmb{v}\right\|_{\pmb{M}}^2 := \pmb{v}^{\mathrm{T}}\pmb{M}\pmb{v}$. Evaluating the functional $J(\pmb{x}_0)$ is very costly because of the nonlinearity of the operators $\mathcal{A}_k$ and $\mathcal{C}_k$. Incre\-mental 4D-Var overcomes this by using consecutive linearised versions of the opera\-tors and minimising a quadratic functional \cite{Courtier1994Incremental}. This approach to minimising~\eqref{4D-Var_Cost} is an inexact Gauss-Newton method \cite{Lawless2005_GaussNewton}. Linearisations of $\mathcal{A}_k$ and $\mathcal{C}_k$ about the model state $\pmb{x}_k$ lead to Jacobians $\pmb{A}_k$ and $\pmb{C}_k$ in the so-called \textit{tangent linear model}:
\begin{align}\label{DynSys4}
    \delta \pmb{x}_{k+1} = \pmb{A}_k \delta \pmb{x}_{k}, \enspace
    \pmb{d}_{k} = \pmb{C}_k \delta \pmb{x}_{k},
\end{align}
where $\pmb{d}_{k} := \pmb{m}_k - \mathcal{C}_k(\pmb{x}_{k})$ for all $k$ and $\delta \pmb{x}_{k}$ denotes a state perturbation. A discussion of how the observation and prior uncertainties influence \eqref{DynSys4} follows in Sect.\@ \ref{BT-4D-Var} The corresponding quadratic cost functional that has to be minimised for $\delta \pmb{x}_0$ is given by
\begin{subequations}\label{4D-Var_InnerCost}
\begin{align}\label{4D-Var_InnerCost_Alone}
\Tilde{J}(\delta \pmb{x}_0) = \frac{1}{2}\left\| \delta \pmb{x}_0 - [ \pmb{x}_0 - \pmb{x}_0^b] \right\|^2_{\pmb{\Gamma}_\mathrm{pr}^{-1}}&+\frac{1}{2}\sum_{k=0}^n \left\| \pmb{C}_k\delta\pmb{x}_k - \pmb{d}_k \right\|^2_{\pmb{R}_k^{-1}}
\end{align}
subject to the forward model dynamics
\begin{align}
\delta \pmb{x}_{k+1} &= \pmb{A}_k\delta \pmb{x}_{k} \enspace \text{for } k = 0,\ldots, n-1.
\end{align}
\end{subequations}
inding the minimum $\delta \pmb{x}_0^*$ of \eqref{4D-Var_InnerCost} is called \textit{inner loop}. The result of the inner loop is the next increment $\delta \pmb{x}_0^{(i)}$ used to update $\pmb{x}_0^{(i+1)} = \pmb{x}_0^{(i)} + \delta \pmb{x}_0^{(i)}$ in the so-called \textit{outer loop}. The trajectory $\pmb{x}^{(i)} = [\pmb{x}_0^{(i)},\pmb{x}_1^{(i)},\ldots,\pmb{x}_n^{(i)}]$ and expected observations are obtained by evaluating the nonlinear model based on the current step's initial condition $\pmb{x}_0^{(i)}$.

Despite the promising idea of obtaining $\pmb{x}_0^*$ incrementally, the computation of \eqref{4D-Var_InnerCost} is still very costly in high-dimensional settings, even in the linear case. For this rea\-son, the system \eqref{DynSys4} should be reduced to obtain lower-dimensional matrices re\-plac\-ing $\pmb{A}_k$, $\pmb{C}_k$ and $\pmb{\Gamma}_\mathrm{pr}$ in \eqref{4D-Var_InnerCost}. Details on how balanced truncation is used to make the min\-imisation of this cost functional  tractable is discussed in Sect.\@ \ref{BT-4D-Var}

\subsection{Model Reduction by Balanced Truncation}\label{BT}
Projection-based methods are popular in model order reduction for dynamical sys\-tems \cite{Benner2015Survey,Benner2017MORBook}. The system-theoretic concept of balanced truncation \cite{Moore1981BT,Mullis1976BT} is one of them. We provide some background from systems theory in Sect.\@ \ref{SystemTheory} before ex\-plaining the proposed method in Sect.\@ \ref{BTAlgo} The theoretical explanations in these sections are based on a book by Antoulas \cite{Antoulas2005Book}. The time-limited approach to balanced truncation  is described in detail in Sect.\@ \ref{TLBT-Defs}

\subsubsection{System-Theoretic Basics}\label{SystemTheory}
Balanced truncation (BT) is usually considered for a linear time invariant (LTI) sys\-tem with continuous ($t \in \mathbb{R}_+$) or discrete ($t \in \mathbb{Z}_+$) state equation as follows:
\begin{subequations}\label{DynSys2}
\begin{align}
\begin{split}\label{DynSys2a}
    \frac{\mathrm{d}\pmb{x}(t) }{\mathrm{d}t}  &= \pmb{Ax}(t) + \pmb{Bu}(t) \enspace \text{ or } \enspace \pmb{x}(t+1) = \pmb{Ax}(t)  + \pmb{Bu}(t),
\end{split}\\
\begin{split}\label{DynSys2b}
    \pmb{y}(t)  &= \pmb{Cx}(t),
\end{split}\\
\begin{split}\label{DynSys2c}
    \pmb{x}(0)  &= \pmb{0}.
\end{split}
\end{align}
\end{subequations}
The output equation and the initial condition are unified for both cases. Here, $\pmb{x},\,\pmb{A}$, and $\pmb{C}$ are as in \eqref{DynSys1}. The input of the system $\pmb{u}(t) \in \mathbb{R}^{d_\mathrm{in}}$ is through the port \linebreak\mbox{$\pmb{B} \in \mathbb{R}^{d \times d_\mathrm{in}}$}. The observer output $\pmb{y}(t) \in \mathbb{R}^{d_\mathrm{out}}$ is without any noise. We, excep\-tional\-ly, use the notation $\pmb{x}(t+1)$ instead of $\pmb{x}_{k+1}$ for discrete time in this subsec\-tion to unify the concept of BT for discrete and continuous time. It is important to notice, however, that the physical meaning of the system matrices $\pmb{A}$ and $\pmb{B}$ differs from continuous time to discrete time: when we discretise the continuous-time sys\-tem $\frac{\mathrm{d}\pmb{x}(t)}{\mathrm{d}t}  = \pmb{A}_{\text{cont}}\pmb{x}(t) + \pmb{B}_{\text{cont}}\pmb{u}(t)$ to $\pmb{x}(t+1) = \pmb{A}_{\text{disc}}\pmb{x}(t)  + \pmb{B}_{\text{disc}}\pmb{u}(t)$, it holds that $\pmb{A}_{\text{disc}} = e^{\pmb{A}_{\text{cont}}}$ and $\pmb{B}_{\text{disc}} = \int_0^1e^{\pmb{A}_{\text{cont}}\tau}\mathrm{d}\tau \pmb{B}_{\text{cont}}$ under the assumption that $\pmb{u}(t)$ is constant on each interval $[k,k+1),\,k\in \mathbb{N}$.

The idea of BT is to truncate the system by removing the states that are difficult to reach and difficult to observe. The concepts of reachability and observability rely on the so-called reachability and observability Gramians, spanning the (sub)spaces of reachable and observable states. The Gramians are functions of time and the \textbf{time-limited reachability Gramian} $\pmb{P}(t) \in \mathbb{R}^{d \times d}$ at time $0 <t < \infty$ is defined by
\begin{subequations}\label{P}
\begin{align}
    \pmb{P}(t) &= \int_0^t e^{\pmb{A}\tau}\pmb{B}\pmb{B}^{\mathrm{T}}e^{\pmb{A}^{\mathrm{T}}\tau}\mathrm{d}\tau \enspace  \text{for }  t \in \mathbb{R}_+ \text{ and }\label{PCont} \\
    \pmb{P}(t) &= \sum_{k=0}^{t-1} \pmb{A}^k\pmb{B}\pmb{B}^{\mathrm{T}}(\pmb{A}^{\mathrm{T}})^k \enspace  \text{for }  t \in \mathbb{Z}_+.\label{PDis}
\end{align}
\end{subequations}
The notion of observability of states at a fixed time $t$ is dual to reachability. The \textbf{time-limited observability Gramian} $\pmb{Q}(t) \in \mathbb{R}^{d \times d}$ at time $0 <t < \infty$ is defined by
\begin{subequations}\label{Q}
\begin{align}
    \pmb{Q}(t) &= \int_0^t e^{\pmb{A}^{\mathrm{T}}\tau}\pmb{C}^{\mathrm{T}}\pmb{C}e^{\pmb{A}\tau}\mathrm{d}\tau \enspace \text{for } t \in \mathbb{R}_+ \text{ and }\label{QCont} \\
    \pmb{Q}(t) &= \sum_{k=0}^{t-1} (\pmb{A}^{\mathrm{T}})^k\pmb{C}^{\mathrm{T}}\pmb{C}\pmb{A}^k \enspace \text{for } t \in \mathbb{Z}_+.\label{QDis}
\end{align}
\end{subequations}
The Gramians are associated with the reachability energy and the observability en\-ergy of a state, respectively:
\begin{align}\label{Energies}
    \left\|\pmb{x}\right\|_{[\pmb{P}(t)]^{-1}}^2 = \pmb{x}^\mathrm{T}[\pmb{P}(t)]^{-1}\pmb{x}, \enspace \enspace\left\|\pmb{x}\right\|_{\pmb{Q}(t)}^2 = \pmb{x}^\mathrm{T}\pmb{Q}(t)\pmb{x},
\end{align}
where $\left\|\pmb{x}\right\|_{[\pmb{P}(t)]^{-1}}^2$ is the minimal energy necessary to reach the state $\pmb{x}$ at time $t$ and $\left\|\pmb{x}\right\|_{\pmb{Q}(t)}^2$ is the maximal energy produced by the output of an observable state $\pmb{x}$ at time $t$. The time-limited Gramians are, thus, functions of time informing us how much energy is necessary to reach or observe a state at a fixed, finite time $t$.

By monitoring the system from initial time to infinity, it is possible to determine how much energy is required to reach or observe a state at \textit{any} point in time. The lim\-its of the finite Gramians \eqref{P} and \eqref{Q} exist only under certain restrictive conditions: the system is required to be \textit{stable}, i.e., the eigenvalues of $\pmb{A}$ lie in the open left half plane for the continuous-time setting or lie inside the unit circle for the discrete-time setting. Then, $e^{\pmb{A}\tau}$ is bounded for $\tau \rightarrow \infty$ and, hence, $\pmb{P}(t)$ and $\pmb{Q}(t)$ are bounded for $t \rightarrow \infty$ and their limits $\pmb{P}_\infty = \lim_{t\rightarrow \infty}  \pmb{P}(t)$ and $\pmb{Q}_\infty  = \lim_{t\rightarrow \infty} \pmb{Q}(t)$ exist; simi\-larly for the discrete-time case. The stable LTI system \eqref{DynSys2} has \textbf{infinite reachability Gramian} $\pmb{P}_\infty \in \mathbb{R}^{d \times d}$ defined by
\begin{subequations}\label{P_inf}
\begin{align}
    \pmb{P}_\infty &= \int_0^\infty e^{\pmb{A}\tau}\pmb{B}\pmb{B}^{\mathrm{T}}e^{\pmb{A}^{\mathrm{T}}\tau}\mathrm{d}\tau\enspace \text{ for } t \in \mathbb{R}_+  \text{ and } \label{P_inf_Cont}\\
     \pmb{P}_\infty &= \sum_{k=0}^\infty \pmb{A}^k\pmb{B}\pmb{B}^{\mathrm{T}}(\pmb{A}^{\mathrm{T}})^k\enspace  \text{ for }  t \in \mathbb{Z}_+.\label{P_inf_Dis}  
\end{align} 
\end{subequations}
Besides, it has \textbf{infinite observability Gramian} $\pmb{Q}_\infty \in \mathbb{R}^{d \times d}$ defined by
\begin{subequations}\label{Q_y}
\begin{align}  
    \pmb{Q}_\infty &= \int_0^\infty
    e^{\pmb{A}^{\mathrm{T}}\tau}\pmb{C}^{\mathrm{T}}\pmb{C}e^{\pmb{A}\tau}\mathrm{d}\tau\enspace  \text{ for }  t \in \mathbb{R}_+  \text{ and } \label{Q_y_Cont} \\
    \pmb{Q}_\infty &= \sum_{k=0}^\infty (\pmb{A}^{\mathrm{T}})^k\pmb{C}^{\mathrm{T}}\pmb{C}\pmb{A}^k\enspace \text{ for }  t \in \mathbb{Z}_+.\label{Q_y_Dis}
\end{align}
\end{subequations}
The advantage of infinite Gramians is that they can be computed efficiently as the solutions to the continuous-time Lyapunov equations  \mbox{$\pmb{A}\pmb{P}_\infty + \pmb{P}_\infty\pmb{A}^{\mathrm{T}} = -\pmb{B}\pmb{B}^{\mathrm{T}}$} and \mbox{$\pmb{A}^{\mathrm{T}}\pmb{Q}_\infty + \pmb{Q}_\infty\pmb{A} = -\pmb{C}^{\mathrm{T}}\pmb{C}$} or to the Lyapunov equations \mbox{$\pmb{A}\pmb{P}_\infty\pmb{A}^{\mathrm{T}} + \pmb{B}\pmb{B}^{\mathrm{T}} = \pmb{P}_\infty$} and \mbox{$\pmb{A}^{\mathrm{T}}\pmb{Q}_\infty\pmb{A}+\pmb{C}^{\mathrm{T}}\pmb{C} = \pmb{Q}_\infty$} in discrete time.
The drawback is that their use limits the application of BT to stable systems. This is our motivation to discuss an approach to BT for unstable systems, which are common in data assimilation applications. 

\subsubsection{Balanced Truncation}\label{BTAlgo}
Projection-based model reduction methods aim to find a suitable low-dimensional subspace of the state space containing the important features. A projection operator projects the state variables and matrices of the system \eqref{DynSys2} onto this low-dimensional subspace. In our case, the low-dimensional subspace consists of the states that are both easily reachable and easily observable, quantified by the energies \eqref{Energies}. The reachability energy is low and the observability energy is high. The directions ob\-tained by BT, therefore, maximise the generalised Rayleigh quotient $\left(\frac{\left\| \pmb{x}\right\|_{\pmb{Q}_\infty}}{\left\|\pmb{x}\right\|_{\pmb{P}_\infty^{-1}}}\right)^2$. The do\-minant eigendirections of the matrix pencil $(\pmb{Q}_\infty,\pmb{P}_\infty^{-1})=:(\pmb{Q},\pmb{P}^{-1})$, i.e., the eigen\-directions $\pmb{v}_k$ satisfying $\pmb{Q}\pmb{v}_k=\delta^2_k\pmb{P}^{-1}\pmb{v}_k$ and associated with the dominant eigen\-values $\delta^2_k$, maximise the generalised Rayleigh quotient. The eigenvalue square roots \mbox{$\delta_1>\delta_2>\ldots>\delta_d$} are called the \textit{Hankel singular values}.

After having determined the dominant directions $\pmb{v}_k$, the goal is to build a \textit{ba\-lancing transformation} $\pmb{T}$ such that $\pmb{T}^{-1}\pmb{P}\pmb{T}^{-\mathrm{T}} = \tilde{\pmb{P}} =\tilde{\pmb{Q}} = \pmb{T}^{\mathrm{T}}\pmb{Q}\pmb{T}$. If $\tilde{\pmb{P}} = \tilde{\pmb{Q}}$, reachability (related to $\tilde{\pmb{P}}$) and observability (related to $\tilde{\pmb{Q}}$) are expressed in a com\-mon basis and can be considered simultaneously. The balancing transformation is used as the projector to reduce the system dimensions. Balanced truncation for stable LTI systems is based on the singular value decomposition and given by Algorithm \ref{AlgoBT}.

\begin{algorithm}[h]
\caption{Square Root Balanced Truncation}\label{AlgoBT}
\hspace*{\algorithmicindent} \textbf{Input}: system matrices $\pmb{A},\pmb{B},\pmb{C}$ of system \eqref{DynSys2}, truncation rank $r\leq d$.\\
\hspace*{\algorithmicindent} \textbf{Output}: reduced system matrices $\Hat{\pmb{A}},\Hat{\pmb{B}},\Hat{\pmb{C}}$, reduced state $\Hat{\pmb{x}}$.\\
\begin{algorithmic}[1]
\State Obtain a square root decomposition (e.g., Cholesky) of the Gramians (\eqref{P_inf}, \eqref{Q_y} or others), i.e., $\pmb{P} = \pmb{RR}^{\mathrm{T}}, \pmb{Q} = \pmb{LL}^{\mathrm{T}}$, e.g., by low-rank Lyapunov solvers \cite{Benner2013NumSolSurvey}.
\State $\pmb{U}\pmb{\Delta}\pmb{Z}^{\mathrm{T}} =\pmb{L}^{\mathrm{T}} \pmb{R}$ (Singular value decomposition).
\State Truncation $\Hat{\pmb{\Delta}} := \mathrm{diag}(\delta_1,\ldots,\delta_r)$, $\Hat{\pmb{U}} := \pmb{U}[:,1:r] $, $\Hat{\pmb{Z}} := \pmb{Z}[:,1:r] $.
\State Balancing transformation $\pmb{T} := \pmb{R}\Hat{\pmb{Z}}\Hat{\pmb{\Delta}}^{-\frac{1}{2}}$, $\pmb{T}^{-1} :=\Hat{\pmb{\Delta}}^{-\frac{1}{2}}\Hat{\pmb{U}}^{\mathrm{T}} \pmb{L}^{\mathrm{T}}$.
\State Reduced system $\Hat{\pmb{x}}:= \pmb{T}^{-1}\pmb{x}$, $\Hat{\pmb{A}} := \pmb{T}^{-1}\pmb{A}\pmb{T}$, $\Hat{\pmb{B}} := \pmb{T}^{-1} \pmb{B}$, $\Hat{\pmb{C}} := \pmb{C}\pmb{T}$.
\end{algorithmic}
\end{algorithm}
\begin{remark}
The matrix $\pmb{\Delta} = \mathrm{diag}(\delta_1,\ldots,\delta_d)$ contains the Hankel singular val\-ues and $\pmb{RZ} =: \pmb{V} = \begin{bmatrix}\pmb{v}_1&\ldots&\pmb{v}_d \end{bmatrix}$ contains the corresponding eigendirections of the pencil $(\pmb{Q},\pmb{P}^{-1})$.
\end{remark}
\begin{remark}
The infinite Gramians $\Hat{\pmb{P}}$, $\Hat{\pmb{Q}} $ of the reduced system
\begin{align}\label{DynSys2_red}
    \frac{\mathrm{d}\Hat{\pmb{x}}}{\mathrm{d}t}  = \Hat{\pmb{A}}\Hat{\pmb{x}} + \Hat{\pmb{B}}\pmb{u}(t), \enspace
  \Hat{\pmb{y}} &= \Hat{\pmb{C}}\Hat{\pmb{x}}
\end{align}
are given by $\Hat{\pmb{P}} = \Hat{\pmb{Q}} = \Hat{\pmb{\Delta}}$ and, thus, the system is balanced.
\end{remark}
\begin{remark}
The balanced truncation Algorithm \ref{AlgoBT} can be applied to LTI systems regardless of details, as long as Gramians are defined. The typical applications are stable systems with infinite Gramians where Gramian square roots are obtained by solving Lyapunov equations. In Sect.\@ \ref{BT for DA} we discuss adaptions of this algorithm to slightly different Gramians. The next subsection introduces time-limited Gramians for unstable systems.
\end{remark}

\subsubsection{Time-limited Balanced Truncation}\label{TLBT-Defs}
Time- and frequency-limited balanced truncation has been developed by Gawronski and Juan \cite{Wodek1990TLGramians}. This section is based on their work. The idea of time-limited balanced truncation is to use the time-limited Gramians \eqref{P} and \eqref{Q} for a time interval of inter\-est instead of the infinite Gramians \eqref{P_inf} and \eqref{Q_y}. In the data assimilation problems, we consider the interval $\mathcal{T} = [0,\ldots, t_e = t_n]$ and use the notation $\pmb{P}_{\mathcal{T}} := \pmb{P}(t_e)$ and $\pmb{Q}_{\mathcal{T}} := \pmb{Q}(t_e)$. The time-limited Gramians are the solutions of (modified) Lyapunov equations:

\begin{theorem}\label{Lyapunov_TLBT}
The time-limited Gramians $\pmb{P}_{\mathcal{T}}$ and $\pmb{Q}_{\mathcal{T}}$ for a continuous LTI system \eqref{DynSys2} with stable $\pmb{A}$ are the unique, positive-semidefinite solutions to the modified conti\-nuous-time Lyapunov equations
\begin{align*}
    \pmb{A}\pmb{P}_{\mathcal{T}} + \pmb{P}_{\mathcal{T}}\pmb{A}^{\mathrm{T}} &= -\pmb{B}\pmb{B}^{\mathrm{T}} + \pmb{B}_{t_e}\pmb{B}_{t_e}^{\mathrm{T}},\enspace \pmb{B}_{t_e} := e^{\pmb{A}t_e}\pmb{B}, \text{ and } \\
   \pmb{A}^{\mathrm{T}}\pmb{Q}_{\mathcal{T}} + \pmb{Q}_{\mathcal{T}}\pmb{A} &= -\pmb{C}^{\mathrm{T}}\pmb{C} + \pmb{C}_{t_e}^{\mathrm{T}}\pmb{C}_{t_e},\enspace \pmb{C}_{t_e} := \pmb{C}e^{\pmb{A}t_e}.
\end{align*}
The time-limited Gramians $\pmb{P}_{\mathcal{T}}$ and $\pmb{Q}_{\mathcal{T}}$ for a discrete system \eqref{DynSys2} with stable $\pmb{A}$ are the unique, positive-semidefinite solutions to the modified discrete-time Lyapunov equations
\begin{align*}
    \pmb{P}_{\mathcal{T}} &= \pmb{A}\pmb{P}_{\mathcal{T}}\pmb{A}^{\mathrm{T}} + \pmb{A}^{t_s}\pmb{B}\pmb{B}^{\mathrm{T}}(\pmb{A}^{\mathrm{T}})^{t_s} - \pmb{A}^{t_e}\pmb{B}\pmb{B}^{\mathrm{T}}(\pmb{A}^{\mathrm{T}})^{t_e}, \text{ and } \\
    \pmb{Q}_{\mathcal{T}} &= \pmb{A}^{\mathrm{T}}\pmb{Q}_{\mathcal{T}}\pmb{A}+(\pmb{A}^{\mathrm{T}})^{t_s}\pmb{C}^{\mathrm{T}}\pmb{C}\pmb{A}^{t_s} - (\pmb{A}^{\mathrm{T}})^{t_e}\pmb{C}^{\mathrm{T}}\pmb{C}\pmb{A}^{t_e}.
\end{align*}
\end{theorem}
The balanced truncation Algorithm \ref{AlgoBT} is applied to the system \eqref{DynSys2} using the Gramians $\pmb{P} = \pmb{P}_{\mathcal{T}}$ and $\pmb{Q} = \pmb{Q}_{\mathcal{T}}$. Their matrix square roots may be obtained as the solutions of modified Lyapunov equations from Theorem \ref{Lyapunov_TLBT}. We call this procedure \textit{time-limited balanced truncation (TLBT)}. Using TLBT is arguably the way to go in the setting of data assimilation, since there are only observations in some fixed time-interval $\mathcal{T} = [0,t_e]$. The end time $t_e$ is often comparatively small and we are not interested in the behaviour at an infinite time horizon.

The assumption that the system \eqref{DynSys2} and hence $\pmb{A}$ is stable is only necessary for the existence of the infinite Gramians $\pmb{P}_\infty$ and $\pmb{Q}_\infty$. The time-limited Gramians $\pmb{P}(t)$ and $\pmb{Q}(t)$ are defined anyways, see \eqref{P} and \eqref{Q}. Hence, we drop the stability assumption for TLBT. Redmann and Kürschner \cite{Redmann2018TLBT} have shown that for continuous systems with unstable $\pmb{A}$ and $\Lambda(\pmb{A}) \cap \Lambda(-\pmb{A}) = \emptyset$ the Lyapunov equations from Theorem \ref{Lyapunov_TLBT} can still be used for the computation of the time-limited Gramians, which is beneficial from a numerical perspective. A drawback is that, even for stable $\pmb{A}$, stability is not preserved for systems \eqref{DynSys2_red}, which were reduced by TLBT, but additional assumptions or modifications are necessary \cite{Wodek1990TLGramians}.

In this section, two data assimilation approaches and the basics of balanced trun\-ca\-tion and its time-limited version have been explained. In the next Sect.\@ \ref{BT for DA}, the details of (TL)BT in these two methods are discussed and compared.

\section{Time-limited Balanced Truncation for Data Assimilation}\label{BT for DA}
The idea of balanced truncation (BT) from Sect.\@ \ref{BT} has been applied to 4D-Var for discrete time \cite{BoessDiss,GreenDiss,Lawless2008Strong4DVarBT} and to linear Gaussian (LG) Bayesian inference for continuous time \cite{Qian2021Balancing}. An adaptation of (time-limited) BT to data assimilation is only possible if the underlying dynamical system can be considered as an LTI system for which Gramians are defined. In the following, we describe how this adaptation is done for the inner loop of 4D-Var (Sect.\@ \ref{BT-4D-Var}) and LG Bayesian inference (Sect.\@ \ref{BT-BLIP}) and introduce time-limited BT for both cases.

\subsection{TLBT within the Inner Loop of Incremental 4D-Var}\label{BT-4D-Var}
BT for the 4D-Var method was first proposed by Lawless et al.\@ \cite{Lawless20084D-Varproposal} and then applied to settings with \cite{GreenDiss} or without \cite{BoessDiss,Lawless2008Strong4DVarBT} model error in the state dynamics \eqref{DynSys3a}. In Sect.\@ \ref{4DVar} we have explained why a reduced order model is needed within the inner loop of incremental 4D-Var.

The system to be truncated is the tangent linear model \eqref{DynSys4}. Let us assume that $\pmb{A}_k = \pmb{A}$ and $\pmb{C}_k = \pmb{C}$ for all $k$. This is an important assumption and restriction, because BT as given in Sect.\@ \ref{BT} is only defined for LTI systems. In \eqref{DynSys4}, the noise has been ignored to give the basic idea of incremental 4D-Var. With noise, the tangent model in the inner loop of 4D-Var taking into account the prior is:
\begin{subequations}\label{DynSys5}
\begin{align}
    \delta \pmb{x}_{0}=\pmb{e}_0 \sim \mathcal{N}(\pmb{0},\pmb{\Gamma}_\mathrm{pr}),\enspace
    \delta \pmb{x}_{k+1} &= \pmb{A} \delta \pmb{x}_{k} \label{DynSys5a},\\
    \pmb{d}_{k} &= \pmb{C} \delta \pmb{x}_{k},\label{DynSys5b}
\end{align}
\end{subequations}
making sure that the covariance of the increment $\delta \pmb{x}_0$ corresponds to the prior covari\-ance \cite{GreenDiss}. Again, we assume the observation error to be time-invariant, i.e., $\pmb{R}_k = \pmb{\Gamma}_{\epsilon}$ for all $k$.
A reduced version of \eqref{DynSys5} is:
\begin{align}\label{DynSys5_red}
\begin{split}
    \delta \Hat{\pmb{x}}_{0}=\Hat{\pmb{e}}_0 \sim \mathcal{N}(\pmb{0},\Hat{\pmb{\Gamma}}_\mathrm{pr}), \enspace
    \delta \Hat{\pmb{x}}_{k+1} &= \Hat{\pmb{A}} \delta \Hat{\pmb{x}}_{k},\\
    \Hat{\pmb{d}}_{k} &= \Hat{\pmb{C}} \delta \Hat{\pmb{x}}_{k},
\end{split}
\end{align}
where the increment $\delta \pmb{x} \in \mathbb{R}^d$ is projected by $\pmb{T}^{-1}\delta \pmb{x}=: \delta \Hat{\pmb{x}} \in \mathbb{R}^r$, with \mbox{$r \ll d$}. The model and observation operators are reduced by $\Hat{\pmb{A}} :=\pmb{T}^{-1}\pmb{A}\pmb{T}$ and \mbox{$\Hat{\pmb{C}}:=\pmb{C}\pmb{T}$}. The projection also affects the state variables, hence $\Hat{\pmb{x}}:=\pmb{T}^{-1}\pmb{x}$. This justifies the projected version of the prior covariance matrix $\Hat{\pmb{\Gamma}}_\mathrm{pr}:=\pmb{T}^{-1}\pmb{\Gamma}_\mathrm{pr}\pmb{T}^{-\mathrm{T}}$; the observa\-tions $\pmb{m}_k$ and the observation error covariance $\pmb{\Gamma}_\epsilon$ remain unchanged. We propose to compute the projection operator $\pmb{T}$ via time-limited balanced truncation.

System \eqref{DynSys5} is almost of type \eqref{DynSys2}, but the input comes with a covariance and thus a different input port at time $k=0$. Green \cite{GreenDiss} proposes to begin the summation for the time-limited Gramians with zero input port at $k=1$. For time step $k=0$ the summand $\pmb{A}^0\pmb{\Gamma}_\mathrm{pr}(\pmb{A}^{\mathrm{T}})^0 = \pmb{\Gamma}_\mathrm{pr}$ is added. With these adjustments, the time-limited Gramians \eqref{PDis} and \eqref{QDis} of \eqref{DynSys5} are:
\begin{subequations}\label{4D-Var-Gramians}
\begin{align}
    \pmb{P}_\mathcal{T}^{\text{4D-Var}} &= \pmb{\Gamma}_\mathrm{pr} + \sum_{k=1}^{t_e-1} \pmb{A}^k\pmb{0}(\pmb{A}^{\mathrm{T}})^k = \pmb{\Gamma}_\mathrm{pr} \, \text{and} \label{4D-Var-P_inf}\\
    \pmb{Q}_\mathcal{T}^{\text{4D-Var}} &= \sum_{k=0}^{t_e-1} (\pmb{A}^{\mathrm{T}})^k\pmb{C}^{\mathrm{T}}\pmb{\Gamma}_{\epsilon}^{-1}\pmb{C}\pmb{A}^k\label{4D-Var-Q_y}.
\end{align}
\end{subequations}
The relevant literature is inconsistent regarding the inclusion of the observation error from \eqref{DynSys3b} in the observability Gramian for the tangent linear model \eqref{DynSys5}. From com\-parison to \eqref{DynSys3b}, it holds that $\pmb{d}_{k} = \pmb{\epsilon}_{k} \sim \mathcal{N}(\pmb{0},\pmb{\Gamma}_{\epsilon})$ for the error $\pmb{d}_{k} := \pmb{m}_k - \mathcal{C}_k(\pmb{x}_{k})$. In the theses of Boess and Green \cite{BoessDiss,GreenDiss} the observability Gramian for \eqref{DynSys5} was given by $\pmb{Q}_\infty^{\text{4D-Var}} = \sum_{k=0}^\infty (\pmb{A}^{\mathrm{T}})^k\pmb{C}^{\mathrm{T}}\pmb{\Gamma}_{\epsilon}\pmb{C}\pmb{A}^k$, which we do not consider correct. Instead, we prefer a time-limited formulation of Lawless et al.\@ \cite{Lawless2008Strong4DVarBT}, which corresponds exactly to \eqref{4D-Var-Q_y}. Bernstein et al.\@ \cite{Bernstein1986Optimal} state that a model reduction by BT with observability Gramian \eqref{4D-Var-Q_y} minimises the distance error between the ouput $\pmb{d}_k$ of the full tangent linear model \eqref{DynSys5} and the output $\Hat{\pmb{d}_k}$ of the reduced model \eqref{DynSys5_red}, weighted by the ob\-servation precision $\pmb{\Gamma}_{\epsilon}^{-1}$, i.e., $\lim_{k\rightarrow \infty} \mathbb{E}\left( (\pmb{d}_k - \Hat{\pmb{d}_k})^{\mathrm{T}}\pmb{\Gamma}_{\epsilon}^{-1}(\pmb{d}_k - \Hat{\pmb{d}_k})\right)$. Weighting by the inverse of the covariance matrix means "penalizing errors in the approximation [...] more strongly in directions of lower [...] variance." \cite[p.\@ 11]{Spantini2015Optimal} This also makes sense here: In directions of low observational covariance, we are fairly confident that our observation is correct. The observation errors  $\pmb{d}_k$ of the full 4D-Var system \eqref{DynSys3} are low and so should their approximations $\Hat{\pmb{d}_k}$. If $|\pmb{d}_k - \Hat{\pmb{d}_k}|$ are large, our approx\-imation is very inaccurate and has to be penalised. In contrast, if we are uncertain about the observations, the errors $\pmb{d}_k$ have larger variance. If the approximations $\Hat{\pmb{d}_k}$ are relatively inaccurate but within the standard deviation, the reduced system \eqref{DynSys5_red} approximates \eqref{DynSys5} reasonably well.

With this explanation, we will from now on use \eqref{4D-Var-Q_y} as observability Gramian for TLBT in incremental 4D-Var. The time-limited observability Gramian for 4D-Var satisfies a modified discrete Lyapunov equation,
\begin{align}\label{4DVarLyap}
\pmb{Q}_{\mathcal{T}}^{\text{4D-Var}} &=\pmb{A}^{\mathrm{T}}\pmb{Q}^{\text{4D-Var}}_{\mathcal{T}}\pmb{A}+\pmb{C}^{\mathrm{T}}\pmb{\Gamma}_{\epsilon}^{-1}\pmb{C}- (\pmb{A}^{\mathrm{T}})^{t_e}\pmb{C}^{\mathrm{T}}\pmb{\Gamma}_{\epsilon}^{-1}\pmb{C}\pmb{A}^{t_e}.
\end{align}
The proof of this equation follows along the same lines as for the deterministic case; see, e.g., \cite{Antoulas2005Book}.

Time-limited balanced truncation is applied to system \eqref{DynSys5} for the inner loop of incremental 4D-Var with Gramians \eqref{4D-Var-Gramians}. The balancing transformation $\pmb{T}$ and its in\-verse $\pmb{T}^{-1}$ from Algorithm \ref{AlgoBT} are the projection matrices to obtain the reduced system \eqref{DynSys5_red}. If the linear tangent model \eqref{DynSys5} does not change during the outer loop iterations, the reduced system \eqref{DynSys5_red} may be computed once with TLBT for the entire incremental 4D-Var algorithm, making this technique particularly effective.

Standard BT in data assimilation is limited to stable systems. We overcome this by considering time-limited BT. Another approach for unstable discrete systems is $\alpha$-bounded balanced truncation, which has been developed for 4D-Var by Boess \cite{Boess20114DVarunstable}. The idea of $\alpha$-bounded BT is to shift the system matrices by a parameter $\alpha$ so that the eigenvalues of $\pmb{A}$ lie inside a disc of radius $\alpha$ around the origin. $\alpha$-bounded BT requires experimentation on the choice of the parameter $\alpha$ \cite{GreenDiss} and is only applicable to discrete systems. We have, therefore, focused on a more classical system-theoretic concept and have introduced time-limited BT. In Sect.\@ \ref{TLBT-Exps} we demonstrate that this approach additionally improves the accuracy for low end times.

\subsection{Balancing Linear Gaussian Bayesian Inference}\label{BT-BLIP}
BT has been extended to continuous linear Gaussian (LG) Bayesian inference by Qian et al.\@ \cite{Qian2021Balancing}. We generalise the concept to arbitrary prior covariances and unstable system matrices $\pmb{A}$. There are two main differences between the LTI system \eqref{DynSys2} and the LG Bayesian inference setting \eqref{DynSys1} that need to be resolved before applying BT:
\begin{enumerate}
      \item In the LG Bayesian inference setting, there is no input function $\pmb{u}(t)$.
      \item Unlike for the LTI system, the observations $\pmb{m}_k$ in Bayesian inference are noisy.
\end{enumerate}
  
\subsubsection{From Compatible to Arbitrary Prior Covariances as Reachability Gramians}\label{BT-DA-pr}
Qian et al.\@ \cite{Qian2021Balancing} have solved the problem of the missing input in Bayesian infer\-ence by introducing the notion of \textit{prior-compatibility}. It is not necessary to know the input port $\pmb{B}$ explicitly but only to ensure that there is a $\pmb{B} \in \mathbb{R}^{d \times d_\mathrm{in}}$ with $d_\mathrm{in} = \mathrm{rank}(\pmb{A}\pmb{\Gamma}_{\mathrm{pr}} + \pmb{\Gamma}_{\mathrm{pr}}\pmb{A}^{\mathrm{T}})$ satisfying $\pmb{A}\pmb{\Gamma}_{\mathrm{pr}} + \pmb{\Gamma}_{\mathrm{pr}}\pmb{A}^{\mathrm{T}} = -\pmb{B}\pmb{B}^{\mathrm{T}}$. A prior covari\-ance matrix $\pmb{\Gamma}_{\mathrm{pr}}$ fulfilling this is called \textit{prior-compatible}. The prior-compatible $\pmb{\Gamma}_{\mathrm{pr}}$ is used as infinite reachability Gramian for the LTI system with given $\pmb{A}$ and possibly unknown $\pmb{B}$. It is also the time-limited reachability Gramian: $\pmb{P}_\mathcal{T}^{\text{LG}} = \pmb{P}_\infty^{\text{LG}} =  \pmb{\Gamma}_{\mathrm{pr}}$.

 Prior-compatibility is not ensured in general. The spectral decomposition of \linebreak\mbox{$\pmb{A}\pmb{\Gamma}_{\mathrm{pr}} + \pmb{\Gamma}_{\mathrm{pr}}\pmb{A}^{\mathrm{T}}$} may be used to manipulate the prior covariance and make it prior-compatible \cite[Sect.\@ 4.1.2]{Qian2021Balancing}. In practice, however, the prior distribution is usually not chosen arbitrarily, but based on expert knowledge or accumulated information. Changing the prior to make it compatible and artificially introducing an input port $\pmb{B}$ may involve a loss of information. We want to justify that any Gaussian prior covari\-ance is a time-limited reachability Gramian.

For the tangent linear model in 4D-Var the time-limited reachability Gramian is given by \eqref{4D-Var-P_inf}. By comparison of the discrete linear tangent model \eqref{DynSys5} and the continuous LG Bayesian inference setting, we can use the continuous version of \eqref{4D-Var-P_inf} to compute the time-limited reachability Gramian for \eqref{DynSys1}, which is
\begin{subequations}\label{BLIP-Gramians}
\begin{align}\label{BLIP-P_inf}
   \pmb{P}_\mathcal{T}^{\text{LG}} = \int_{0}^{t_e} e^{\pmb{A}\tau}\pmb{BB}^\mathrm{T}e^{\pmb{A}^{\mathrm{T}}\tau}\mathrm{d}\tau + \pmb{\Gamma}_\mathrm{pr} =  \int_{0}^{t_e} e^{\pmb{A}\tau} \pmb{00}^{\mathrm{T}}  e^{\pmb{A}^{\mathrm{T}}\tau}\mathrm{d}\tau + \pmb{\Gamma}_\mathrm{pr} =  \pmb{\Gamma}_{\mathrm{pr}}.
\end{align}
The above comparison with established results from 4D-Var (see Sect.\@ \ref{BT-4D-Var}) over\-comes the need for prior-compatibility for $\pmb{\Gamma}_{\mathrm{pr}}$ and makes the application of BT to LG Bayesian inference more general. A numerical example to emphasise that the original prior should be used instead of a modified compatible one whenever possible is given in Sect.\@ \ref{ExperimentDescription}

This approach is also well-founded in literature: time-limited reachability Grami\-ans are the covariance matrices of the LTI system state at a fixed time \cite{Melina2018SDE} and the infinite Gramian is, thus, the total state covariance. The state equation \eqref{BLIP_forward} evolves linearly without model error. Hence, the state evolution is deterministic and the full state covariance is introduced by the prior covariance as reachability Gramian. The POD-Gramian of an input-free system is constructed similarly \cite{Antoulas2005Book}.

\begin{remark}
    The unknown initial condition to be inferred is not zero. It has a cen\-tered Gaussian prior and a Gaussian posterior distribution. BT is defined for a zero initial condition. It may therefore be questioned whether this method is well suited for application to data assimilation problems. Other systems theory techniques, however, require an input matrix; see, e.g., \cite{Antoulas2005Book,Benner2015Survey,Benner2017MORBook}. This input would have to be artificially con\-structed by modifying prior knowledge, which is not necessary for BT. Therefore, we consider balanced truncation to be the best method so far.
\end{remark}

\subsubsection{A Noisy Time-Limited Observability Gramian}
Regarding the pollution of the measurement $\pmb{m}_k$ by noise $\pmb{\epsilon}_k \sim \mathcal{N}(\pmb{0},\pmb{\Gamma}_\epsilon)$ in the LG Bayesian inference setting \eqref{DynSys1} compared to the LTI system \eqref{DynSys2}, the following is proposed \cite[p.\@ 28]{Qian2021Balancing}: Recall that $\left\|\pmb{x}\right\|^2_{\pmb{Q}}=\left\|\pmb{y}\right\|^2_{L^2(\mathbb{R})}$ defines the maximal energy produced by an observable state $\pmb{x}$. It is the squared distance of the output signal $\pmb{y}$ to $\pmb{0}$ in the $L^2(\mathbb{R})$-norm. For noisy observations, instead of $\left\|\pmb{\cdot}\right\|^2_{L^2(\mathbb{R})}$, the Mahalanobis distance \cite{Mahalanobis1936} from the equilibrium state $\pmb{0}$ to the conditional distribution of a single measurement $\pmb{m}_k|(\pmb{x}_0,\pmb{t}_k) \sim \mathcal{N}(\pmb{C}e^{\pmb{A}t_k}\pmb{x}_0, \pmb{\Gamma}_\epsilon)$ is used, i.e., $$d\left(\pmb{0},\mathcal{N}(\pmb{C}e^{\pmb{A}t_k}\pmb{x}_0, \pmb{\Gamma}_\epsilon)\right) = \left\|\pmb{C}e^{\pmb{A}t_k}\pmb{x}_0\right\|_{\pmb{\Gamma}_\epsilon^{-1}}.$$
We use the same reasoning for the definition of the \textbf{noisy time-limited observability Gramian} $\pmb{Q}_\mathcal{T}^{\text{LG}}\in \mathbb{R}^{d\times d}$ for linear Gaussian Bayesian inference as:
\begin{align}\label{BLIP-Q_m}
    \pmb{Q}_\mathcal{T}^{\text{LG}} = \int_0^{t_e} e^{\pmb{A}^{\mathrm{T}}\tau}\pmb{C}^{\mathrm{T}}\pmb{\Gamma}_\epsilon^{-1}\pmb{C}e^{\pmb{A}\tau}\mathrm{d}\tau.
\end{align}
\end{subequations}
It is the unique, positive-definite solution to the dual modified Lyapunov equation
\begin{align}\label{BLIP-Q_m-Lyap}
    \pmb{A}^{\mathrm{T}}\pmb{Q}_\mathcal{T}^{\text{LG}} + \pmb{Q}_\mathcal{T}^{\text{LG}}\pmb{A} = -\pmb{C}^{\mathrm{T}}\pmb{\Gamma}_\epsilon^{-1}\pmb{C} + \pmb{C}_{t_e}^{\mathrm{T}}\pmb{\Gamma}_\epsilon^{-1}\pmb{C}_{t_e},\enspace \pmb{C}_{t_e} := \pmb{C}e^{\pmb{A}t_e},
\end{align}
similar to the discrete Lyapunov equation \eqref{4DVarLyap} for 4D-Var. Remark that $\pmb{Q}_\mathcal{T}^{\text{4D-Var}}$ is also a \textit{noisy} (discrete) time-limited observability Gramian. The Fisher matrix \linebreak\mbox{$\pmb{H} = \sum_{k=1}^n e^{\pmb{A}^{\mathrm{T}}t_k}\pmb{C}^{\mathrm{T}}\pmb{\Gamma}_{\pmb{\epsilon}}^{-1}\pmb{C}e^{\pmb{A}t_k}$} is a discrete summation approach to this same \linebreak Gramian and could be used as such \cite{Qian2021Balancing}. This resembles an optimal low-rank method for LG Bayesian inference \cite{Spantini2015Optimal} not based on systems theory.

For (TL)BT, Algorithm \ref{AlgoBT} is applied to the LG Bayesian inference problem by using the Gramians \eqref{BLIP-Gramians} in Step 1. This leads to the reduced system
\begin{align}\label{DynSyn1-red}
\begin{split}
    \frac{d\Hat{\pmb{x}}}{dt}  &= \Hat{\pmb{A}}\Hat{\pmb{x}}, \text{ especially }\Hat{\pmb{x}}(0) = \pmb{T}^{-1}\pmb{x}_0,\\
  \pmb{m}_k &\approx \Hat{\pmb{C}}\Hat{\pmb{x}} + \pmb{\epsilon}_k,\enspace \pmb{\epsilon}_k \sim \mathcal{N}(\pmb{0},\pmb{\Gamma}_\epsilon) \text{ as before },
\end{split}
\end{align}
with reduced forward map $\pmb{G}_{\mathrm{(TL)BT}}$ and approximate Fisher information $\pmb{H}_{\mathrm{(TL)BT}}$:
\begin{align*}
\pmb{G}_{\mathrm{(TL)BT}} = \begin{bmatrix}\Hat{\pmb{C}}e^{\Hat{\pmb{A}}t_1}\pmb{T}^{-1}\\\vdots\\\Hat{\pmb{C}}e^{\Hat{\pmb{A}}t_n}\pmb{T}^{-1}\end{bmatrix} =  \begin{bmatrix}\Hat{\pmb{C}}e^{\Hat{\pmb{A}}t_1}\\\vdots\\\Hat{\pmb{C}}e^{\Hat{\pmb{A}}t_n}\end{bmatrix}\pmb{T}^{-1},\enspace
\pmb{H}_{\mathrm{(TL)BT}} = \pmb{G}_{\mathrm{(TL)BT}}^{\mathrm{T}}\pmb{\Gamma}_{\mathrm{obs}}^{-1}\pmb{G}_{\mathrm{(TL)BT}},
\end{align*}
and gives the posterior mean and covariance approximations:
    
\begin{equation}
    \begin{aligned}\label{PostStats-red}
        \pmb{\mu}_{\mathrm{pos,(TL)BT}} &= \pmb{\Gamma}_{\mathrm{pos,(TL)BT}}\pmb{G}_{\mathrm{(TL)BT}}^{\mathrm{T}}\pmb{\Gamma}_{\mathrm{obs}}^{-1}\pmb{m},\\ \pmb{\Gamma}_{\mathrm{pos,(TL)BT}} &= (\pmb{H}_{\mathrm{(TL)BT}} + \pmb{\Gamma}_{\mathrm{pr}}^{-1})^{-1}.
    \end{aligned}
\end{equation}
The notation $(\pmb{\mu}_{\mathrm{pos,TLBT}},\pmb{\Gamma}_{\mathrm{pos,TLBT}})$ denotes a reduced posterior obtained by TLBT and $(\pmb{\mu}_{\mathrm{pos,BT}},\pmb{\Gamma}_{\mathrm{pos,BT}})$ a reduced posterior obtained by BT with infinite Gramians.

\begin{remark}
    For initial state $\pmb{x}(0)=\pmb{0}$ there exists an error bound for the output of the BT-reduced system \eqref{DynSyn1-red} compared to the full system \cite{Antoulas2005Book}. For TLBT, error bounds have been derived for continuous time \cite{Redmann2020TLBT,Redmann2018TLBT} and discrete time \cite{PontesDuff2021DiscreteTLBT}. These results are hard to adapt to Bayesian inference, since we infer an unknown non-zero initial state $\pmb{x}(0) = \pmb{x}_0$.
\end{remark}

\begin{remark}
  BT for inhomogeneous, i.e., non-zero, initial conditions has been inves\-tigated by Heinkenschloss et al. proposing to augment the input matrix of the LTI system \eqref{DynSys2} by the initial condition before BT is applied \cite{Heink2011inhom}. This technique seems to be less useful for model reduction with changing or unknown initial conditions. Beattie et al.\@ use splitting: the output is approximated by the simultaneous reduction of two systems, one depending on the initial state $\pmb{x}_0 \neq \pmb{0}$ and one with homogeneous initial condition \cite{Beattie2016inhom}. Regarding the application to Bayesian inference, we have a dif\-ferent objective: not the output approximation, but the reduced posterior statistics.
\end{remark}

\subsubsection{Computational Aspects of TLBT in Data Assimilation}\label{TLBT-Comp}
The use of the time-limited noisy Gramians \eqref{4D-Var-Gramians} and \eqref{BLIP-Gramians} in Algorithm \ref{AlgoBT} is straight\-forward if in Step 1 there are (low-rank) factors for $\pmb{P}_{\mathcal{T}} = \pmb{RR}^{\mathrm{T}}$ and \mbox{$\pmb{Q}_{\mathcal{T}} = \pmb{LL}^{\mathrm{T}}$}. In LG Bayesian inference and 4D-Var, the time-limited reachability Gramian is equal to the prior covariance. Therefore, no Lyapunov equation needs to be solved. The time-limited noisy observability Gramian computation requires the solution of the modified Lyapunov equation \eqref{4DVarLyap} or \eqref{BLIP-Q_m-Lyap}. The main computational obstacle in solv\-ing the discrete or continuous modified Lyapunov equations is the computation of the matrix power $\pmb{A}^{t_e}$ or the matrix exponential $e^{\pmb{A}t_e}$ for the right-hand sides. There is no direct square root approach. But direct approaches quickly reach their limit for high-dimensional problems. Slower iterative algorithms, like (block-)rational Krylov algorithms, may still compute a result. Shift selection is a typical bottleneck in these algorithms and they take longer than solving the Lyapunov equations for the infinite Gramians in low-dimensional settings; see, e.g., \cite{Kuerschner2018TLBT,PontesDuff2021DiscreteTLBT}. Time-limited Gramians, however, allow the application of the BT Algorithm \ref{AlgoBT} to unstable systems, where Lyapunov equation solvers cannot be used.

\begin{table}[h]
\footnotesize
    \begin{tabular*}{\textwidth}{l|c|c}
    System   &   TL reachability Gramian  &  TL observability Gramian\\
    \midrule
    \multirow{2.5}{2.2cm}{Standard TLBT \eqref{DynSys2}} &  $\pmb{P}_{\mathcal{T}} = \int_0^{t_e} e^{\pmb{A}\tau}\pmb{B}\pmb{B}^{\mathrm{T}}e^{\pmb{A}^{\mathrm{T}}\tau}\mathrm{d}\tau$ & $\pmb{Q}_{\mathcal{T}} = \int_0^{t_e} e^{\pmb{A}^{\mathrm{T}}\tau}\pmb{C}^{\mathrm{T}}\pmb{C}e^{\pmb{A}\tau}\mathrm{d}\tau$\\
    \cmidrule(rl){2-3}
     & $\pmb{P}_{\mathcal{T}} = \sum_{k=0}^{t_e-1} \pmb{A}^k\pmb{B}\pmb{B}^{\mathrm{T}}(\pmb{A}^{\mathrm{T}})^k$ & $\pmb{Q}_{\mathcal{T}} = \sum_{k=0}^{t_e-1} (\pmb{A}^{\mathrm{T}})^k\pmb{C}^{\mathrm{T}}\pmb{C}\pmb{A}^k$\\ 
    \midrule
    LG Infer. \eqref{DynSys1} & $\pmb{P}_{\mathcal{T}}^{\text{LG}} =  \pmb{\Gamma}_{\mathrm{pr}}$ &$\pmb{Q}_{\mathcal{T}}^{\text{LG}} = \int_0^{t_e} e^{\pmb{A}^{\mathrm{T}}\tau}\pmb{C}^{\mathrm{T}}\pmb{\Gamma}_\epsilon^{-1}\pmb{C}e^{\pmb{A}\tau}\mathrm{d}\tau$ \\ 
    \midrule
   4D-Var \eqref{DynSys5} & $\pmb{P}_{\mathcal{T}}^{\text{4D-Var}} =\pmb{\Gamma}_{\mathrm{pr}}$ &$\pmb{Q}_{\mathcal{T}}^{\text{4D-Var}} = \sum_{k=0}^{t_e-1} (\pmb{A}^{\mathrm{T}})^k\pmb{C}^{\mathrm{T}}\pmb{\Gamma}_{\epsilon}^{-1}\pmb{C}\pmb{A}^k$ \\
    \end{tabular*}
    \caption{Time-limited (TL) Gramians in data assimilation problems}
    \label{TableGramians}
\end{table}
In this section, time-limited balanced truncation has been introduced for the lin\-ear tangent model of 4D-Var and for LG Bayesian inference. Table \ref{TableGramians} summarises the resulting Gramians for TLBT in data assimilation. By comparing the two systems in question, we have generalised the application of BT in linear Gaussian Bayesian inference to arbitrary prior covariances. We have provided a unifying definition for the observability Gramian in incremental 4D-Var. Our main contribution is that our approach can handle unstable system matrices $\pmb{A}$, which are common in data assim\-ilation problems. In Sect.\@ \ref{TLBT-Exps} we will see that TLBT further outperforms the standard BT approach for short end-times.

\section{Numerical Experiments for TLBT in Data Assimilation}\label{TLBT-Exps}
This section is dedicated to demonstrating the performance of TLBT in linear Gaus\-sian Bayesian inference with numerical experiments. The MATLAB code reproduc\-ing our results is available at \url{https://github.com/joskoUP/TLBTforDA}.

\subsection{Comparison of Compatible and Non-Compatible Prior Covariance}\label{ExperimentDescription}
In the previous section we have argued that it is reasonable to set $\pmb{P}_\mathcal{T}^{\text{LG}} = \pmb{\Gamma}_{\mathrm{pr}}$ regard\-less of the prior-compatibility. The behaviour of non-compatible versus compatible priors in BT for LG Bayesian inference will be demonstrated by a numerical ex\-periment. The reduced posterior statistics are given by \eqref{PostStats-red}. We include the optimal low-rank approach (OLR) to the posterior statistics \eqref{PostStats} by Spantini et al.\@ \cite{Spantini2015Optimal}:

Consider the linear Gaussian Bayesian inference problem $\pmb{m} = \pmb{G}\pmb{x}+ \pmb{\epsilon}$ with for\-ward matrix \mbox{$\pmb{G} \in \mathbb{R}^{nd_\mathrm{out} \times d}$}, Fischer matrix $\pmb{H} \in \mathbb{R}^{d \times d}$ and a vector of measurements \linebreak\mbox{$\pmb{m} = [\pmb{m}_1^\mathrm{T},\ldots,\pmb{m}_n^\mathrm{T}]^\mathrm{T}$}. Prior and observation errors are assumed to be independent, centered Gaussians of the right dimensions with covariance matrices $\pmb{\Gamma}_{\mathrm{pr}}$ and $\pmb{\Gamma}_{\epsilon}$, re\-spectively. The posterior statistics of $\mathcal{N}(\pmb{\mu}_{\mathrm{pos}},\pmb{\Gamma}_{\mathrm{pos}})$ are approximated by $\Hat{\pmb{\mu}}_{\mathrm{pos}}$ and $\Hat{\pmb{\Gamma}}_{\mathrm{pos}}$ as follows: 

Let $(\tau_i^2,\pmb{v}_i)_{i=1}^d$ be the $d$ generalised eigenvalue-eigenvector pairs  of the matrix pencil $(\pmb{H},\pmb{\Gamma}_{\mathrm{pr}}^{-1})$ and let $(\pmb{w}_i)_{i=1}^d$ be the $d$ generalised eigenvectors of the matrix pencil $(\pmb{G}\pmb{\Gamma}_\mathrm{pr}\pmb{G}^\mathrm{T},\pmb{\Gamma}_\mathrm{obs})$.  Both sequences of generalised eigenvectors are assumed to be associated with a non-increasing sequence of eigenvalues and the $\pmb{w}_i$ are nor\-malised w.r.t.\@ $\pmb{\Gamma}_\mathrm{obs}$. The optimal posterior covariance approximation $\Hat{\pmb{\Gamma}}_{\mathrm{pos}}$ of $\pmb{\Gamma}_{\mathrm{pos}}$ in the Förstner-metric (out of class $\mathcal{M}_r:=\{\pmb{\Gamma}_{\mathrm{pr}}-\pmb{KK}^{\mathrm{T}} \succ 0\,|\, \mathrm{rank}(\pmb{K}) \leq r\}$) is
\begin{subequations}\label{Spantini}
    \begin{align}
        \Hat{\pmb{\Gamma}}_{\mathrm{pos}} = \pmb{\Gamma}_{\mathrm{pr}} - \sum_{i=1}^r \frac{\tau_i^2}{1+ \tau_i^2}\pmb{v}_i\pmb{v}_i^{\mathrm{T}}.
    \end{align}
The optimal mean approximation $\Hat{\pmb{\mu}}_\mathrm{pos}$ of $\pmb{\mu}_\mathrm{pos}$ in the $\left\|\cdot\right\|_{\pmb{\Gamma}_{\mathrm{pos}}^{-1}}$-norm Bayes risk (out of class $\mathcal{V}_r^{LR}:=\{\pmb{v} = \pmb{Nm}\,|\,\pmb{N}\in \mathbb{R}^{d\times nd_\mathrm{out}},\,\mathrm{rank}(\pmb{N})\leq r\}$) is
    \begin{align}
        \Hat{\pmb{\mu}}_\mathrm{pos} = \sum_{i=1}^r \frac{\tau_i}{1+ \tau_i^2}\pmb{v}_i\pmb{w}_i^{\mathrm{T}}\pmb{m}.
    \end{align}
\end{subequations}
See the paper by Spantini et al.\@ \cite{Spantini2015Optimal} for more details and explanations of the measures.
\begin{remark}
\begin{enumerate}[label=(\alph*)] 
    \item For semi-positive definite matrices $\pmb{A}$ and $\pmb{B}$ and $(\sigma_i)_i$ the gener\-alised eigenvalues of the matrix pencil $(\pmb{A},\pmb{B})$, the Förstner-metric $d_{\mathcal{F}} (\pmb{A},\pmb{B})$ is defined as $d_{\mathcal{F}}^2 (\pmb{A},\pmb{B}):= \mathrm{trace}\big[ \ln^2(\pmb{A}^{-\frac{1}{2}}\pmb{B}\pmb{A}^{-\frac{1}{2}}) \big] = \sum_i \ln^2(\sigma_i)$ \cite{Foerstner2003}.
    \item The $\left\|\cdot\right\|_{\pmb{\Gamma}_{\mathrm{pos}}^{-1}}$-norm Bayes risk $R(\delta(\pmb{m}),\pmb{x})$ of the measurement-based estimator $\delta(\pmb{m})$ of parameter $\pmb{x}$ unknown is defined as $R(\delta(\pmb{m}),\pmb{x}) := \mathbb{E}(\left\|\pmb{x} - \delta(\pmb{m})\right\|^2_{{\pmb{\Gamma}_{\mathrm{pos}}^{-1}}})$. The expectation is taken over the joint distribution of $\pmb{x}$ and $\pmb{m}$ \cite{LehmCase1998PointEst}.
\end{enumerate}
\end{remark} 
For the numerical experiments in this and the next subsection we have reproduced the benchmark examples by Qian et al.\@ \cite[cf.\@ p.\@ 22]{Qian2021Balancing}.\footnote{Code from \url{https://github.com/elizqian/balancing-bayesian-inference}.} In the first presented ex\-periment, we have modified the code for the additional computations with the non-compatible prior. We always consider LTI systems starting at $t = 0$ and running up to different end times $t_e$. Measurements are taken at discrete equidistant times $t_i = ih$, \mbox{$(i = 1,\ldots,n)$,} $t_n = t_e$. The goal is to infer the initial condition $\pmb{x}_0$ by us\-ing reduced models of different ranks $r$. For our experiments, a true initial condition is drawn from $\mathcal{N}(\pmb{0},\pmb{\Gamma}_{\mathrm{pr}})$ and measurements are generated from evolving the under\-lying dynamical system exactly and adding $\mathcal{N}(\pmb{0},\pmb{\Gamma}_{\epsilon})$-Gaussian measurement noise. We consider a given prior non-compatible with the system matrix $\pmb{A}$. Using the code from \cite[Appendix D]{Qian2021Balancing}, we create a modified compatible prior.

The reduced posterior statistics \eqref{PostStats-red} for this problem are computed using BT with infinite Gramians for both the given non-compatible (NC) and the modified compati\-ble (C) prior. We compare them with the full posterior statistics \eqref{PostStats} using the Förstner metric (for the covariance) and the $\left\|\cdot\right\|_{\pmb{\Gamma}_{\mathrm{pos}}^{-1}}$-norm Bayes risk (for the mean) derived from the true non-compatible prior (-NC) and derived with the modified compati\-ble prior (-C). OLR for both priors and both measures is also considered. For the Bayes risk computation, we run 100 independent experiments to obtain 100 different posterior mean estimates and compute the empirical Bayes risk among them.

We study the heat equation\footnote{LTI system matrices and documentation at \url{http://slicot.org/20-site/126-benchmark-examples-for-model-reduction}.} with $d = 200$, $d_\mathrm{out}=1$ and $\pmb{\Gamma}_{\epsilon}= \sigma_{\mathrm{obs}}^2$, \linebreak\mbox{$\sigma_{\mathrm{obs}} = 0.008$.} The non-compatible prior covariance $\pmb{\Gamma}_{\mathrm{pr}} = \pmb{R}\pmb{R}^{\mathrm{T}}$ with $\pmb{R} = \linebreak\mathrm{diag}(\mathrm{ones}(1,d)) + \mathrm{diag}(0.5\cdot \mathrm{randn}(1,d-1),1) + \mathrm{diag}(0.25\cdot \mathrm{randn}(1,d-2),2)$ takes into account that neighbouring points might influence one another with decreas\-ing intensity.
\begin{figure}[h]
    
\definecolor{mycolor1}{rgb}{0.85000,0.32500,0.09800}%
\definecolor{mycolor2}{rgb}{1.00000,0.00000,1.00000}%
\definecolor{mycolor3}{rgb}{0.00000,1.00000,1.00000}%
\begin{tikzpicture}
\footnotesize
\begin{axis}[%
width=0.36\textwidth,
height=1.1in,
at={(3.327in,0.481in)},
scale only axis,
xmin=0,
xmax=20,
xlabel style={font=\color{white!15!black}},
xlabel={$r$},
ymode=log,
ymin=1e-15,
ymax=1000,
yminorticks=true,
axis background/.style={fill=white},
title={$\pmb{\Gamma}_\mathrm{pos}$ F\"orstner error},
xmajorgrids,
ymajorgrids,
yminorgrids,
legend style={at={(0,0)}, anchor=south west, legend cell align=left, align=left, fill=none, draw=none}
]

\addplot [color=mycolor2, only marks, mark=o, mark options={solid, mycolor2}]
  table[row sep=crcr]{%
1	666.931471240608\\
2	467.7625930678\\
3	336.095857397918\\
4	265.787930673734\\
5	219.70976822708\\
6	167.50566898012\\
7	135.00484076232\\
8	113.560968778948\\
9	98.6137657540848\\
10	88.5617629592481\\
11	83.7072941787174\\
12	78.9917200228025\\
13	77.6776889225972\\
14	77.9991236135382\\
15	78.2790906833449\\
16	78.4058644883823\\
17	78.4061311580947\\
18	78.3899767193732\\
19	78.4005399365897\\
20	78.4081332678264\\
};
\addlegendentry{BT C-NC}

\addplot [color=mycolor2, only marks, mark=+, mark options={solid, mycolor2}]
  table[row sep=crcr]{%
1	577.530052527222\\
2	378.713947035445\\
3	247.173839393909\\
4	176.924347118895\\
5	130.928662957322\\
6	78.9075934440012\\
7	46.5803328683585\\
8	26.0182374681428\\
9	12.5075349903341\\
10	4.62805683910138\\
11	1.45800220290173\\
12	0.551651387882239\\
13	0.160145323469347\\
14	0.0489301183997832\\
15	0.0100232175494767\\
16	0.00122087933195202\\
17	0.000217050574051631\\
18	8.11843676554812e-05\\
19	2.76346581705317e-05\\
20	1.04657289751246e-05\\
};
\addlegendentry{BT C-C}

\addplot [color=mycolor1, only marks, mark=o, mark options={solid, mycolor1}]
  table[row sep=crcr]{%
1	313.522319146233\\
2	239.865556601968\\
3	180.391685345656\\
4	129.618028808565\\
5	79.2330415416473\\
6	45.4786917210636\\
7	25.2790754560075\\
8	11.9772614397212\\
9	4.49926848176295\\
10	1.22968290676921\\
11	0.242554490288472\\
12	0.0843076888773824\\
13	0.0324575557582086\\
14	0.0113098281961773\\
15	0.00486088943782721\\
16	0.000542596833976891\\
17	0.000218085663868187\\
18	0.000102471142666466\\
19	1.34377722825786e-05\\
20	3.20024631722945e-06\\
};
\addlegendentry{BT NC-NC}

\addplot [color=mycolor3, only marks, mark=+, mark options={solid, mycolor3}]
  table[row sep=crcr]{%
1	574.346887892082\\
2	372.680497483394\\
3	233.813720636155\\
4	139.683947579281\\
5	78.9691631509426\\
6	39.0099185331021\\
7	16.8342381244096\\
8	5.29911945079166\\
9	0.913622421583285\\
10	0.186872039635177\\
11	0.0164168407249408\\
12	0.000942177653464374\\
13	4.59624416783822e-05\\
14	1.66418283846892e-06\\
15	2.78308647862536e-08\\
16	1.20798221582152e-09\\
17	7.48701159510419e-11\\
18	2.9564205968524e-12\\
19	1.11379072807458e-13\\
20	1.62774859872765e-15\\
};
\addlegendentry{OLR C-C}

\addplot [color=blue, only marks, mark=o, mark options={solid, blue}]
  table[row sep=crcr]{%
1	309.848216975214\\
2	211.096114045289\\
3	135.484145598038\\
4	78.8065149648828\\
5	39.1404566158341\\
6	16.5542545699996\\
7	5.45122777682404\\
8	1.06711966952522\\
9	0.12799860766017\\
10	0.00624183484040842\\
11	0.000397260584754347\\
12	4.02115556077107e-05\\
13	4.63670534413682e-06\\
14	2.65951931950324e-07\\
15	6.92447772011345e-09\\
16	9.36031410791486e-11\\
17	6.6534584741307e-12\\
18	2.93389522351751e-13\\
19	1.56737583923124e-14\\
};
\addlegendentry{OLR NC-NC}
\end{axis}

\begin{axis}[%
width=0.36\textwidth,
height=1.1in,
at={(0.758in,0.481in)},
scale only axis,
xmin=0,
xmax=20,
xlabel style={font=\color{white!15!black}},
xlabel={$r$},
ymode=log,
ymin=1e-8,
ymax=1,
yminorticks=true,
axis background/.style={fill=white},
title={$\pmb{\mu}_\mathrm{pos}$  Bayes risk},
xmajorgrids,
ymajorgrids,
legend style={at={(0,0)}, anchor=south west, legend cell align=left, align=left, fill=none, draw=none}
]
\addplot [color=mycolor2, only marks, mark=o, mark options={solid, mycolor2}]
  table[row sep=crcr]{%
1	0.386052896123708\\
2	0.244495175308872\\
3	0.152177827031674\\
4	0.0750293662336379\\
5	0.0428586090712055\\
6	0.0405906836906833\\
7	0.0167587390225417\\
8	0.0119473749596965\\
9	0.00941472043402385\\
10	0.00862713339665551\\
11	0.00946636680306332\\
12	0.0131620320070414\\
13	0.0138698896462454\\
14	0.0139751200633639\\
15	0.0137434560207199\\
16	0.0136721449786195\\
17	0.0136485579908086\\
18	0.0136436347666283\\
19	0.0136396111235283\\
20	0.0136366051818882\\
};
\addlegendentry{BT C-NC}

\addplot [color=mycolor2, only marks, mark=+, mark options={solid, mycolor2}]
  table[row sep=crcr]{%
2	0.912063928046916\\
3	0.347467160962135\\
4	0.148590428768008\\
5	0.0752582526521725\\
6	0.0570796181810423\\
7	0.025640844781125\\
8	0.0165874863614475\\
9	0.0087665424685526\\
10	0.00494490050464045\\
11	0.00329695618189865\\
12	0.00222245850851765\\
13	0.00121241838927408\\
14	0.000675757723236743\\
15	0.00040588586838798\\
16	0.000193162027242659\\
17	0.000106996365493973\\
18	4.05582172555051e-05\\
19	2.18851085066846e-05\\
20	1.26463891459059e-05\\
};
\addlegendentry{BT C-C}

\addplot [color=mycolor1, only marks, mark=o, mark options={solid, mycolor1}]
  table[row sep=crcr]{%
1	0.476025575778451\\
2	0.286703444811405\\
3	0.14068586871895\\
4	0.0761948997800316\\
5	0.0689664084461538\\
6	0.0294809600127341\\
7	0.0157804027331334\\
8	0.00902454319351572\\
9	0.00736590823118702\\
10	0.00601898278210265\\
11	0.00239583739872427\\
12	0.00111599459514829\\
13	0.000717892675318768\\
14	0.000461329390517625\\
15	0.000276817868458326\\
16	0.000132523489921974\\
17	0.0001005136830781\\
18	5.6880589763983e-05\\
19	1.68924649715461e-05\\
20	9.30292065481706e-06\\
};
\addlegendentry{BT NC-NC}

\addplot [color=mycolor3, only marks, mark=+, mark options={solid, mycolor3}]
  table[row sep=crcr]{%
2	0.879308949077778\\
3	0.327538025849765\\
4	0.11681348508525\\
5	0.0577788250007668\\
6	0.0280703411398903\\
7	0.0140319113216798\\
8	0.0071819735125198\\
9	0.00340017700885336\\
10	0.00190639845847058\\
11	0.000949595651242844\\
12	0.000466600161138759\\
13	0.000203067760226427\\
14	0.000100085078629704\\
15	2.9093873952132e-05\\
16	1.45176727428239e-05\\
17	7.40079610193969e-06\\
18	3.13106731730403e-06\\
19	1.47139011203676e-06\\
20	4.78916359771495e-07\\
};
\addlegendentry{OLR C-C}

\addplot [color=blue, only marks, mark=o, mark options={solid, blue}]
  table[row sep=crcr]{%
1	0.37103825028395\\
2	0.195238507100506\\
3	0.0974496355085576\\
4	0.055830361523656\\
5	0.0286266046336019\\
6	0.0136504202268316\\
7	0.00718061394908279\\
8	0.00355206774915068\\
9	0.0016917464411878\\
10	0.000686414002294835\\
11	0.000395159930956076\\
12	0.000209349124212759\\
13	0.000111123778019847\\
14	5.66629975469468e-05\\
15	2.10595900105173e-05\\
16	7.85155551636284e-06\\
17	3.87561355444844e-06\\
18	1.82422249483893e-06\\
19	7.60624372899446e-07\\
20	3.24618885246607e-07\\
};
\addlegendentry{OLR NC-NC}
\end{axis}
\end{tikzpicture}%    
       
    \vspace{-3ex}
    \caption{Comparison of a given non-compatible prior (NC) and a modified compatible prior (C) in balanced truncation (BT) \cite{Qian2021Balancing} and the optimal low-rank approach (OLR) \cite{Spantini2015Optimal} of LG Bayesian inference  for the heat equation. The first letters denote the prior used for the approx.\@ posterior computation and the second the one used for the reference. Measurements are spaced $h=0.005$ apart inside $\mathcal{T} = [0,10]$}\label{IncompPr}
 \end{figure}
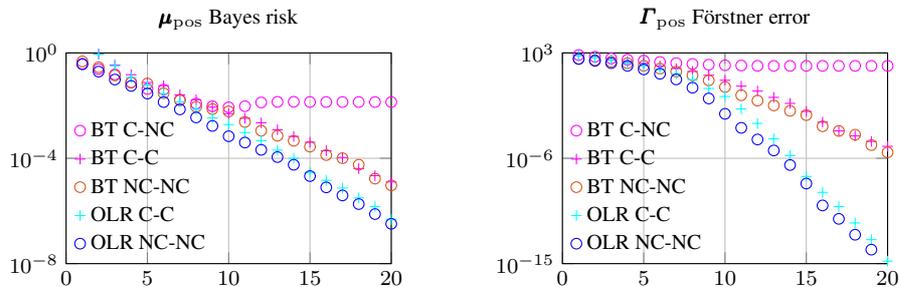

The results of the experiments are depicted in Fig.\@ \ref{IncompPr} Similarly to the results of Qian et al.\@ \cite{Qian2021Balancing}, OLR (in blue/cyan) is superior to BT (pink/orange) for the same prior. BT with the true NC prior (orange circles) performs almost as good as OLR, but the modified compatible prior gives bad approximations in the original (NC) mea\-sures (pink circles). Considering these approximations in the measures with regard to the modified compatible prior (C), they are almost as good as the original ones (pink and cyan crosses). It is thus not useful to enforce prior-compatibility for non-compatible priors containing valuable knowledge as the modification changes the re\-sults significantly. Compatible priors, however, perform well in the measures derived from them. Building compatible priors may therefore be a way of constructing a prior that fits the system dynamics when good prior knowledge is not available otherwise.

Non-compatible priors are particularly important when dealing with unstable sys\-tem matrices $\pmb{A}$. Time-limited BT is especially suitable for such systems in data as\-similation, as will be demonstrated in the next experiments.

\subsection{Comparison of TLBT and Standard BT with Infinite Gramians}

We apply TLBT to benchmark examples and compare the approximate posterior statistics with full balancing (BT) and with the optimal low-rank approach (OLR). For the experiments in this section we consider the same setup as before with added code for the time-limited Gramian computation and with the following change: $\pmb{\Gamma}_{\mathrm{pr}}$ is obtained as the solution of $\pmb{A}\pmb{\Gamma}_{\mathrm{pr}} + \pmb{\Gamma}_{\mathrm{pr}}\pmb{A}^{\mathrm{T}} = -\pmb{B}\pmb{B}^{\mathrm{T}}$ for $\pmb{B}$ according to the experi\-ment; see \cite{Qian2021Balancing}. The TLBT-reduced posterior statistics \eqref{PostStats-red} are computed as described in Sect.\@  \ref{BT-BLIP}

The true posterior statistics \eqref{PostStats} are compared in the $\left\|\cdot\right\|_{\pmb{\Gamma}_{\mathrm{pos}}^{-1}}$-norm Bayes risk and the Förstner-metric with $(\pmb{\mu}_{\mathrm{pos,TLBT}},\pmb{\Gamma}_{\mathrm{pos,TLBT}})$, $(\pmb{\mu}_{\mathrm{pos,BT}},\pmb{\Gamma}_{\mathrm{pos,BT}})$, and \linebreak$(\Hat{\pmb{\mu}}_{\mathrm{pos}},\Hat{\pmb{\Gamma}}_{\mathrm{pos}})$ from OLR. The time-limited Gramians for TLBT are computed by the rational Krylov method by Kürschner \cite[p.\@ 1830]{Kuerschner2018TLBT}.\footnote{Code from \url{https://zenodo.org/record/7366026}.} We have to choose parameters: the approximation tolerance (tol) for the low-rank solution of $f(\pmb{A})\pmb{B}$ and the max\-imum number of rational Krylov iterations (maxit). The Lyapunov equation for the recent low-rank approximation is solved at every $s$-th step. We studied both the heat equation and the ISS1R module.\footnote{LTI system matrices and documentation at \url{http://slicot.org/20-site/126-benchmark-examples-for-model-reduction}} All parameters are given in Table \ref{ParametersStable}.
\begin{table}[h]
    \footnotesize
    %\centering
    \begin{tabular*}{\textwidth}{c|c|c|c|c|c|c|c}
    Example    &   $d$    &  $d_\mathrm{out}$ & $\pmb{B}$ & $\pmb{\Gamma}_{\epsilon}$&   $s$    &  tol & maxit\\
    \midrule
    heat equation &   200    &  1 & $\pmb{I}_d$& $\begin{bmatrix}0.008^2\end{bmatrix}$&   5    &  $10^{-8}$ & 100\\ \midrule
    \multirow{2}{2cm}[5pt]{ISS1R module}    &   270    &  3 & see\footnotemark[\value{footnote}] & $\begin{bmatrix} 0.0025^2 & 0 & 0 \\ 0 & 0.0005^2 & 0 \\ 0 & 0 & 0.0005^2\end{bmatrix}$
    &   3    &  $10^{-4}$ & 100\\
    \end{tabular*}
    \caption{Parameters for the LTI examples}
    \label{ParametersStable}
\end{table}
\vspace{-10pt}
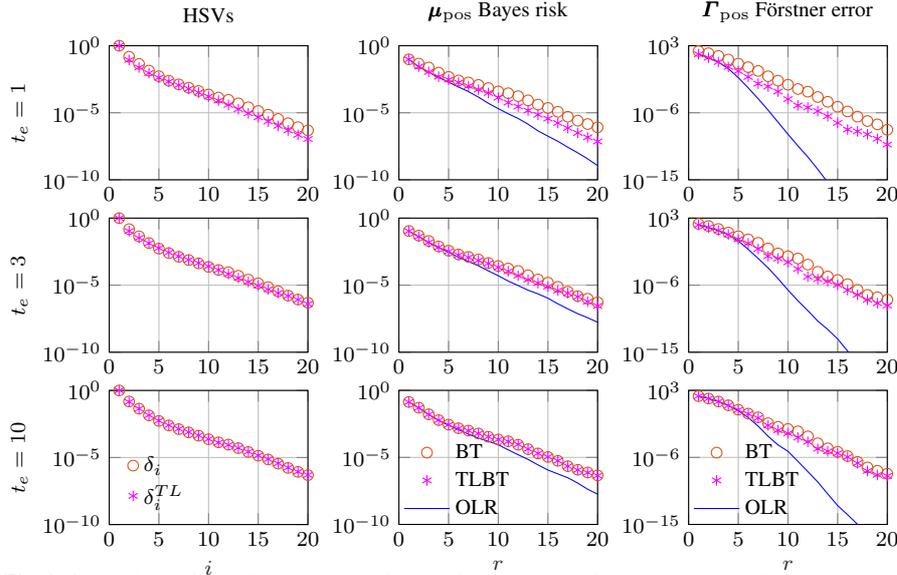
\begin{figure}[h]
\definecolor{mycolor1}{rgb}{0.85000,0.32500,0.09800}%
\definecolor{mycolor2}{rgb}{1.00000,0.00000,1.00000}%
\footnotesize
\begin{tikzpicture}
\begin{axis}[%
width=0.22\textwidth,
height=0.7in,
at={(3.0in,1.8in)},
scale only axis,
xmin=0,
xmax=20,
ymode=log,
ymin=1e-15,
ymax=1000,
yminorticks=true,
axis background/.style={fill=white},
title={$\pmb{\Gamma}_\mathrm{pos}$ F\"orstner error},
xmajorgrids,
ymajorgrids,
yminorgrids
]
\addplot [color=blue, forget plot]
  table[row sep=crcr]{%
1	71.3390351301456\\
2	24.193807454403\\
3	6.03278737434943\\
4	0.820603216439802\\
5	0.0522059658541471\\
6	0.00219217587306647\\
7	7.67476051482327e-05\\
8	2.23296787342153e-06\\
9	5.48303575748843e-08\\
10	1.28259456533597e-09\\
11	3.59527032771935e-11\\
12	1.1438085397146e-12\\
13	2.54667098471413e-14\\
14	4.47909799026103e-16\\
};
\addplot [color=mycolor1, only marks, mark=o, mark options={solid, mycolor1}, forget plot]
  table[row sep=crcr]{%
1	186.854081442156\\
2	99.135585673592\\
3	42.5747792644323\\
4	13.7979461858512\\
5	3.94119783802568\\
6	1.30698056527088\\
7	0.412005681485942\\
8	0.079818626617594\\
9	0.0202573922691126\\
10	0.00579167335281043\\
11	0.00154219563493581\\
12	0.000677743316107106\\
13	0.000129712489063385\\
14	2.95039232400386e-05\\
15	6.42988726631863e-06\\
16	1.49595036429894e-06\\
17	3.04011126215215e-07\\
18	7.78012199867314e-08\\
19	2.76522161798888e-08\\
20	5.22501299014836e-09\\
};
\addplot [color=mycolor2, only marks, mark=asterisk, mark options={solid, mycolor2}, forget plot]
  table[row sep=crcr]{%
1	71.355178376769\\
2	24.3237891525821\\
3	6.4319407112531\\
4	1.41088456841395\\
5	0.391952274479213\\
6	0.0643148825081726\\
7	0.00806179617113574\\
8	0.00339075611020797\\
9	0.000703847593835036\\
10	6.78980222978106e-05\\
11	9.62473435420144e-06\\
12	5.09117561664089e-06\\
13	1.66992566933846e-06\\
14	3.37981220863982e-07\\
15	4.29519509296396e-08\\
16	5.62051342102754e-09\\
17	3.19424960127228e-09\\
18	1.27194308995928e-09\\
19	3.2954553625901e-10\\
20	5.32584020465283e-11\\
};
\end{axis}

\begin{axis}[%
width=0.22\textwidth,
height=0.7in,
at={(1.5in,1.8in)},
scale only axis,
xmin=0,
xmax=20,
ymode=log,
ymin=1e-10,
ymax=1,
yminorticks=true,
axis background/.style={fill=white},
title={$\pmb{\mu}_\mathrm{pos}$  Bayes risk},
xmajorgrids,
ymajorgrids,
yminorgrids
]
\addplot [color=blue, forget plot]
  table[row sep=crcr]{%
1	0.0943931961239697\\
2	0.0273576263609955\\
3	0.0106427310124122\\
4	0.00429770613129688\\
5	0.00192382745647378\\
6	0.000833903920842938\\
7	0.00037934917676656\\
8	0.000146840960493442\\
9	5.09194220649405e-05\\
10	2.03153907722044e-05\\
11	8.4041278785443e-06\\
12	3.7784855559297e-06\\
13	1.39643346292502e-06\\
14	4.45150358294673e-07\\
15	1.78833362691724e-07\\
16	7.45490674816704e-08\\
17	2.67672329948514e-08\\
18	1.02809912769194e-08\\
19	3.63210849842855e-09\\
20	1.13529078078739e-09\\
};
\addplot [color=mycolor1, only marks, mark=o, mark options={solid, mycolor1}, forget plot]
  table[row sep=crcr]{%
1	0.0966512277055191\\
2	0.045935015807419\\
3	0.0237448952688923\\
4	0.00980758935753899\\
5	0.00472770389702511\\
6	0.00254036779520402\\
7	0.00181302935392434\\
8	0.00119059314474634\\
9	0.000695303308849299\\
10	0.000393197698927355\\
11	0.000285306588214236\\
12	0.000172595853436784\\
13	8.39800567631325e-05\\
14	4.06191220498421e-05\\
15	2.28833419946868e-05\\
16	1.18923546470845e-05\\
17	5.06198818052099e-06\\
18	2.57040748918333e-06\\
19	1.43881367526663e-06\\
20	8.23333638875797e-07\\
};
\addplot [color=mycolor2, only marks, mark=asterisk, mark options={solid, mycolor2}, forget plot]
  table[row sep=crcr]{%
1	0.0949923169101168\\
2	0.028135395992146\\
3	0.0111345055026215\\
4	0.00490194357962899\\
5	0.00293641415925385\\
6	0.00180159177682354\\
7	0.000985153982963567\\
8	0.000531766525634756\\
9	0.000301879822450586\\
10	0.000133361478449036\\
11	5.73197627405863e-05\\
12	2.92289275675077e-05\\
13	1.46974010141538e-05\\
14	7.17267455038397e-06\\
15	3.65865615194535e-06\\
16	1.71761208807352e-06\\
17	7.26464867087072e-07\\
18	3.12338187758093e-07\\
19	1.46634947143948e-07\\
20	7.21524660653475e-08\\
};
\end{axis}

\begin{axis}[%
width=0.22\textwidth,
height=0.7in,
at={(0in,1.8in)},
scale only axis,
xmin=0,
xmax=20,
ymode=log,
ymin=1e-10,
ymax=1,
yminorticks=true,
ylabel={$t_e = 1$},
axis background/.style={fill=white},
title={HSVs},
xmajorgrids,
ymajorgrids,
yminorgrids
]
\addplot [color=mycolor1, only marks, mark=o, mark options={solid, mycolor1}, forget plot]
  table[row sep=crcr]{%
1	1\\
2	0.150331219638593\\
3	0.043886338760099\\
4	0.0137883217677465\\
5	0.00529942457881452\\
6	0.00244198749320122\\
7	0.00130235174226813\\
8	0.000745200389818012\\
9	0.000420949266947824\\
10	0.000229461519118154\\
11	0.000135384731799028\\
12	9.18230079694952e-05\\
13	5.01187792545181e-05\\
14	2.57408732339691e-05\\
15	1.36156984514401e-05\\
16	7.12475377460242e-06\\
17	3.48999279531035e-06\\
18	1.66077638754813e-06\\
19	8.47642108280366e-07\\
20	4.85835017438117e-07\\
};
\addplot [color=mycolor2, only marks, mark=asterisk, mark options={solid, mycolor2}, forget plot]
  table[row sep=crcr]{%
1	1\\
2	0.08614805143944\\
3	0.0239920034790834\\
4	0.00912818034487591\\
5	0.00413235802893404\\
6	0.00210524835567188\\
7	0.00113874548698259\\
8	0.000613977138870929\\
9	0.000321847418279581\\
10	0.000163484098049217\\
11	8.06068267865109e-05\\
12	3.88582262058016e-05\\
13	1.8583298119062e-05\\
14	9.03517755913981e-06\\
15	4.53626056965151e-06\\
16	2.2519296188239e-06\\
17	1.05760776392985e-06\\
18	4.81796991328669e-07\\
19	2.2525584560893e-07\\
20	1.10178447521185e-07\\
};
\end{axis}

\begin{axis}[%
width=0.22\textwidth,
height=0.7in,
at={(3.0in,0.9in)},
scale only axis,
xmin=0,
xmax=20,
ymode=log,
ymin=1e-15,
ymax=1000,
yminorticks=true,
axis background/.style={fill=white},
xmajorgrids,
ymajorgrids,
yminorgrids
]
\addplot [color=blue, forget plot]
  table[row sep=crcr]{%
1	121.775543204867\\
2	55.0557271476907\\
3	21.4171020944836\\
4	5.85292020929216\\
5	0.862884482459881\\
6	0.0633636608760323\\
7	0.00385893272967482\\
8	0.000210339541435252\\
9	7.54007355189047e-06\\
10	2.20776195952006e-07\\
11	8.06820062158242e-09\\
12	3.7675164532454e-10\\
13	1.49145153206531e-11\\
14	1.14963574943883e-12\\
15	5.94166790262807e-14\\
16	1.45745888009345e-15\\
17	4.768996446612e-17\\
};
\addplot [color=mycolor1, only marks, mark=o, mark options={solid, mycolor1}, forget plot]
  table[row sep=crcr]{%
1	170.309940855823\\
2	86.6542832961732\\
3	34.7448213976786\\
4	11.6024183636244\\
5	3.66455564278963\\
6	1.41139881747096\\
7	0.486247068843813\\
8	0.0848970741590919\\
9	0.0182755641684162\\
10	0.00765864961873555\\
11	0.00383361642422767\\
12	0.000579625958517783\\
13	6.21129988364023e-05\\
14	1.48957166850403e-05\\
15	5.01362459779284e-06\\
16	1.37420355933554e-06\\
17	2.37344313100101e-07\\
18	4.85105419278002e-08\\
19	2.36473879058408e-08\\
20	1.12516724094985e-08\\
};
\addplot [color=mycolor2, only marks, mark=asterisk, mark options={solid, mycolor2}, forget plot]
  table[row sep=crcr]{%
1	121.782360727186\\
2	55.1147961541594\\
3	21.6040551991048\\
4	6.27863859612649\\
5	1.44945322532077\\
6	0.408281621423373\\
7	0.0725581224778624\\
8	0.00782578093453402\\
9	0.00374463543004205\\
10	0.00117420188911344\\
11	0.000151390833522782\\
12	1.34887022849503e-05\\
13	6.41617074904302e-06\\
14	3.3242537298787e-06\\
15	1.09316925864117e-06\\
16	2.1685237353877e-07\\
17	5.15114591268212e-08\\
18	1.07555507663642e-08\\
19	3.73679038566161e-09\\
20	1.8252251269389e-09\\
};
\end{axis}

\begin{axis}[%
width=0.22\textwidth,
height=0.7in,
at={(1.5in,0.9in)},
scale only axis,
xmin=0,
xmax=20,
ymode=log,
ymin=1e-10,
ymax=1,
yminorticks=true,
axis background/.style={fill=white},
xmajorgrids,
ymajorgrids,
yminorgrids
]
\addplot [color=blue, forget plot]
  table[row sep=crcr]{%
1	0.106861545281555\\
2	0.0410493427590098\\
3	0.014990507584281\\
4	0.00628556977368596\\
5	0.00290236286708666\\
6	0.00122043179940669\\
8	0.000278365330022842\\
9	0.000117085327829076\\
10	5.04172020963357e-05\\
11	2.03687770456111e-05\\
12	9.33104594451442e-06\\
13	4.2126900722661e-06\\
14	2.24279321868153e-06\\
15	1.05841828200802e-06\\
16	4.05431945231233e-07\\
17	1.69747500206792e-07\\
18	7.40242120083948e-08\\
19	3.69344095386255e-08\\
20	1.71036555097503e-08\\
};
\addplot [color=mycolor1, only marks, mark=o, mark options={solid, mycolor1}, forget plot]
  table[row sep=crcr]{%
1	0.110062238595938\\
2	0.0494085863292449\\
3	0.0198160188829906\\
4	0.00771218557086045\\
5	0.00372218253748056\\
6	0.00195433744174925\\
7	0.00128355161611682\\
8	0.000819140008873372\\
9	0.000450795987489573\\
10	0.000267989018600115\\
11	0.000175471177547768\\
12	0.000114458596763519\\
13	5.47312793761353e-05\\
14	2.62925737910976e-05\\
15	1.54389511129132e-05\\
16	7.50346487707095e-06\\
17	3.27515589969835e-06\\
18	1.51552492954049e-06\\
19	9.70456803425139e-07\\
20	5.24784775115568e-07\\
};
\addplot [color=mycolor2, only marks, mark=asterisk, mark options={solid, mycolor2}, forget plot]
  table[row sep=crcr]{%
1	0.107670801696062\\
2	0.0418733102121204\\
3	0.0157615152372142\\
4	0.0067294767082141\\
5	0.00328782312533324\\
6	0.00178170091052667\\
7	0.00117096087384464\\
8	0.00069947397255926\\
9	0.000416623089103389\\
10	0.000212148756743552\\
11	0.000107843210853043\\
12	5.17970920492179e-05\\
13	2.53268821483102e-05\\
14	1.43933706653652e-05\\
15	7.73708832707168e-06\\
16	3.9511806522223e-06\\
17	2.55243490755041e-06\\
18	1.37115138857658e-06\\
19	5.85647640906532e-07\\
20	2.84361480609214e-07\\
};
\end{axis}

\begin{axis}[%
width=0.22\textwidth,
height=0.7in,
at={(0in,0.9in)},
scale only axis,
xmin=0,
xmax=20,
ymode=log,
ymin=1e-10,
ymax=1,
yminorticks=true,
ylabel={$t_e = 3$},
axis background/.style={fill=white},
xmajorgrids,
ymajorgrids,
yminorgrids
]
\addplot [color=mycolor1, only marks, mark=o, mark options={solid, mycolor1}, forget plot]
  table[row sep=crcr]{%
1	1\\
2	0.150331219638593\\
3	0.043886338760099\\
4	0.0137883217677465\\
5	0.00529942457881452\\
6	0.00244198749320122\\
7	0.00130235174226813\\
8	0.000745200389818012\\
9	0.000420949266947824\\
10	0.000229461519118154\\
11	0.000135384731799028\\
12	9.18230079694952e-05\\
13	5.01187792545181e-05\\
14	2.57408732339691e-05\\
15	1.36156984514401e-05\\
16	7.12475377460242e-06\\
17	3.48999279531035e-06\\
18	1.66077638754813e-06\\
19	8.47642108280366e-07\\
20	4.85835017438117e-07\\
};
\addplot [color=mycolor2, only marks, mark=asterisk, mark options={solid, mycolor2}, forget plot]
  table[row sep=crcr]{%
1	1\\
2	0.107966921567328\\
3	0.0333109969924275\\
4	0.0134218298031141\\
5	0.00581355892169305\\
6	0.00273500141365004\\
7	0.00141482660227987\\
8	0.000785049635144582\\
9	0.000449286401656364\\
10	0.000253760401066892\\
11	0.000134734137071951\\
12	6.79007048941384e-05\\
13	3.37898020194257e-05\\
14	1.74178634007339e-05\\
15	9.3513177477563e-06\\
16	4.86006443823761e-06\\
17	2.59202070228117e-06\\
18	1.57700433563707e-06\\
19	8.24707669486755e-07\\
20	3.94681789007365e-07\\
};
\end{axis}

\begin{axis}[%
width=0.22\textwidth,
height=0.7in,
at={(3.0in,0in)},
scale only axis,
xmin=0,
xmax=20,
xlabel style={font=\color{white!15!black}},
xlabel={$r$},
ymode=log,
ymin=1e-15,
ymax=1000,
yminorticks=true,
axis background/.style={fill=white},
xmajorgrids,
ymajorgrids,
yminorgrids,
legend style={at={(-0.02,-0.03)}, anchor=south west, legend cell align=left, align=left, fill=none, draw=none}
]
\addplot [color=mycolor1, only marks, mark=o, mark options={solid, mycolor1}]
  table[row sep=crcr]{%
1	167.198905948678\\
2	80.5763929753339\\
3	31.725923656825\\
4	9.84960586426437\\
5	2.52293198599916\\
6	0.757817401850446\\
7	0.240572827979463\\
8	0.0388205808875085\\
9	0.0124327425384091\\
10	0.00651317597971353\\
11	0.00317895024452543\\
12	0.000666593912281839\\
13	7.85344840218442e-05\\
14	1.65691927552523e-05\\
15	5.62428106667746e-06\\
16	1.49399901341981e-06\\
17	2.5354347222349e-07\\
18	4.47894471788652e-08\\
19	1.61731673577308e-08\\
20	6.46415769645558e-09\\
};
\addlegendentry{BT}

\addplot [color=mycolor2, only marks, mark=asterisk, mark options={solid, mycolor2}]
  table[row sep=crcr]{%
1	165.895804662848\\
2	79.3320980136589\\
3	31.1858959782963\\
4	9.6016712837333\\
5	2.40683620268758\\
6	0.691116688194646\\
7	0.195533766022464\\
8	0.0201571612419516\\
9	0.00437260001267611\\
10	0.00171969080017369\\
11	0.000465484042978828\\
12	0.000123534916596634\\
13	1.46054426613516e-05\\
14	6.64511046841659e-06\\
15	3.61878142213242e-06\\
16	1.11725716073418e-06\\
17	1.68915533569086e-07\\
18	1.54281668204468e-08\\
19	4.45770735025759e-09\\
20	3.21780848151651e-09\\
};
\addlegendentry{TLBT}

\addplot [color=blue]
  table[row sep=crcr]{%
1	165.892202759201\\
2	79.3084311790955\\
3	31.0873991400204\\
4	9.2810719122696\\
5	1.83599304916015\\
6	0.224720565460664\\
7	0.0157867309845407\\
8	0.000691884082661778\\
9	4.72816987007033e-05\\
10	6.56533995200079e-06\\
11	2.95787088487459e-07\\
12	1.25454284542125e-08\\
13	5.00398037128057e-10\\
14	1.24778425767184e-11\\
15	3.37304256015435e-13\\
16	2.01817878120012e-14\\
17	1.02464449117827e-15\\
18	4.56128911768411e-17\\
};
\addlegendentry{OLR}
\end{axis}

\begin{axis}[%
width=0.22\textwidth,
height=0.7in,
at={(1.5in,0in)},
scale only axis,
xmin=0,
xmax=20,
xlabel style={font=\color{white!15!black}},
xlabel={$r$},
ymode=log,
ymin=1e-10,
ymax=1,
yminorticks=true,
axis background/.style={fill=white},
xmajorgrids,
ymajorgrids,
yminorgrids,
legend style={at={(-0.02,-0.03)}, anchor=south west, legend cell align=left, align=left, fill=none, draw=none}
]
\addplot [color=mycolor1, only marks, mark=o, mark options={solid, mycolor1}]
  table[row sep=crcr]{%
1	0.135139677707156\\
2	0.050896290877553\\
3	0.0164527595693386\\
4	0.00605314726616183\\
5	0.00277222162622554\\
6	0.00150293048893865\\
7	0.000996319492715896\\
8	0.000603418629631032\\
9	0.000344350377828952\\
10	0.000216185842809087\\
11	0.000159408527488635\\
12	8.84669617923214e-05\\
13	4.09821856806164e-05\\
14	1.96934054775137e-05\\
15	1.07646082092564e-05\\
16	5.55983262823331e-06\\
17	2.20565982924107e-06\\
18	1.19505891093042e-06\\
19	7.49610427944776e-07\\
20	4.5338167735635e-07\\
};
\addlegendentry{BT}

\addplot [color=mycolor2, only marks, mark=asterisk, mark options={solid, mycolor2}]
  table[row sep=crcr]{%
1	0.134338790991544\\
2	0.0499076235284125\\
3	0.0158852423066967\\
4	0.00595653203760135\\
5	0.00271995253303391\\
6	0.00148691823987047\\
7	0.000989547038079015\\
8	0.000605872997863217\\
9	0.000348639476731171\\
10	0.000214778050646729\\
11	0.000152320897268593\\
12	9.2923929747812e-05\\
13	4.22610567994314e-05\\
14	1.99859678829509e-05\\
15	1.09285839725094e-05\\
16	5.61999813899811e-06\\
17	2.20095784568856e-06\\
18	1.15518878613583e-06\\
19	7.22466129286033e-07\\
20	4.3218029878589e-07\\
};
\addlegendentry{TLBT}

\addplot [color=blue]
  table[row sep=crcr]{%
1	0.133877720948195\\
2	0.0494383486916278\\
3	0.0157022238063612\\
4	0.00572535258604346\\
5	0.00250293932882902\\
6	0.00116135335737303\\
7	0.000543992120812074\\
8	0.000269907662059427\\
9	0.000158547873744303\\
10	8.78129186185085e-05\\
11	3.65718034164186e-05\\
12	1.5825292487056e-05\\
13	7.07299751640309e-06\\
14	2.74530584320495e-06\\
15	1.14428909058214e-06\\
16	5.83917165824592e-07\\
17	3.0369011392814e-07\\
18	1.2135298945768e-07\\
19	4.19550086496325e-08\\
20	1.81255971135595e-08\\
};
\addlegendentry{OLR}

\end{axis}

\begin{axis}[%
width=0.22\textwidth,
height=0.7in,
at={(0in,0in)},
scale only axis,
xmin=0,
xmax=20,
xlabel style={font=\color{white!15!black}},
xlabel={$i$},
ymode=log,
ymin=1e-10,
ymax=1,
yminorticks=true,
ylabel={$t_e = 10$},
axis background/.style={fill=white},
xmajorgrids,
ymajorgrids,
yminorgrids,
legend style={at={(0.05,0.05)}, anchor=south west, legend cell align=left, align=left, fill=none, draw=none}
]
\addplot [color=mycolor1, only marks, mark=o, mark options={solid, mycolor1}]
  table[row sep=crcr]{%
1	1\\
2	0.150331219638593\\
3	0.043886338760099\\
4	0.0137883217677465\\
5	0.00529942457881452\\
6	0.00244198749320122\\
7	0.00130235174226813\\
8	0.000745200389818012\\
9	0.000420949266947824\\
10	0.000229461519118154\\
11	0.000135384731799028\\
12	9.18230079694952e-05\\
13	5.01187792545181e-05\\
14	2.57408732339691e-05\\
15	1.36156984514401e-05\\
16	7.12475377460242e-06\\
17	3.48999279531035e-06\\
18	1.66077638754813e-06\\
19	8.47642108280366e-07\\
20	4.85835017438117e-07\\
};
\addlegendentry{$\delta_i$}

\addplot [color=mycolor2, only marks, mark=asterisk, mark options={solid, mycolor2}]
  table[row sep=crcr]{%
1	1\\
2	0.141315877028199\\
3	0.043657013426522\\
4	0.0142661485186383\\
5	0.00553642924697053\\
6	0.00254994746918486\\
7	0.00135676481223405\\
8	0.000778850629933069\\
9	0.000441327488407186\\
10	0.000238335074976809\\
11	0.000129929689186024\\
12	8.57785773505909e-05\\
13	5.14475661172218e-05\\
14	2.6792313108923e-05\\
15	1.41839384915042e-05\\
16	7.47205597854258e-06\\
17	3.67980341405498e-06\\
18	1.7463284659253e-06\\
19	8.7431768299994e-07\\
20	4.94247527376444e-07\\
};
\addlegendentry{$\delta^{TL}_i$}
\end{axis}
\end{tikzpicture}%
    \vspace{-3ex}
     \caption{Comparison of time-limited balanced truncation (TLBT) with standard balanced truncation (BT) \cite{Qian2021Balancing} and the optimal low-rank approach (OLR) \cite{Spantini2015Optimal} of LG Bayesian inference for the heat equation model. Measurements are spaced $h=0.005$ apart inside $\mathcal{T} = [0,t_e]$ for three different end times $t_e = 1,3,10$. In the left panel, the normalised square roots $\delta_i$ and $\delta^{TL}_i$ of the generalised eigenvalues of the matrix pencils $(\pmb{Q}_{\infty}^{\text{LG}}, \pmb{\Gamma}_{\mathrm{pr}}^{-1})$ and  $(\pmb{Q}_{\mathcal{T}}^{\text{LG}}, \pmb{\Gamma}_{\mathrm{pr}}^{-1})$ are plotted, corresponding to the Hankel singular values}\label{3x3 heat}
\end{figure}
\begin{figure}[h]
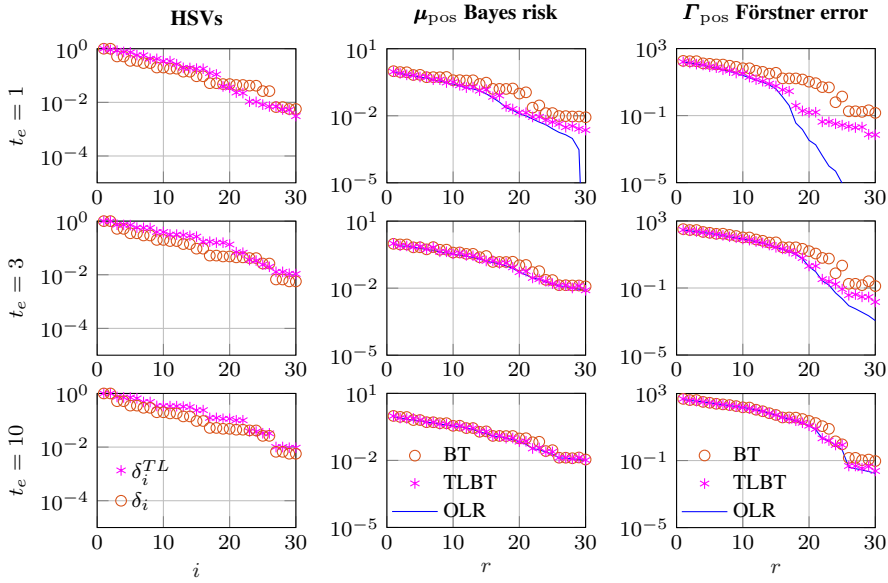

    \centering
\definecolor{mycolor1}{rgb}{0.85000,0.32500,0.09800}%
\definecolor{mycolor2}{rgb}{1.00000,0.00000,1.00000}%
%   
    \vspace{-0.5ex}
    \caption{Quantities plotted as in Figure \ref{3x3 heat}; here for the ISS1R model. Measurements are equispaced with $h=0.1$ inside $\mathcal{T} = [0,t_e]$ for three different end times $t_e = 1,3,10$}\label{3x3 ISS}
 \end{figure}
 
The results depicted in Fig.\@ \ref{3x3 heat} for the heat equation and in Fig.\@ \ref{3x3 ISS} for the ISS model are as expected: for low end times, especially $t_e = 1$, TLBT performs sig\-nificantly better than BT with infinite Gramians. This is caused by the faster decay of the Hankel singular values \cite{Kuerschner2018TLBT}. With increasing end time $t_e$ the TLBT approxi\-mation results approach the BT approximation results. For the heat equation, OLR outperforms TLBT, whereas for the ISS model, optimality can be attained by TLBT, especially for lower reduction ranks~$r$. For low end times, TLBT is able to reach the optimal OLR approximation even when standard BT cannot. This is a huge advan\-tage of using time-limited Gramians. Especially in the ISS1R model for $t_e = 1$ this behaviour can be observed.

OLR is superior to TLBT regarding approximation quality, but not regarding speed. The main advantage of (TL)BT over OLR is that $\pmb{G}$ and $\pmb{H}$ are computed in reduced versions of lower dimensionality. They are obtained by forwarding a reduced rather than a full dynamical system \eqref{DynSys1}. Computations based on a reduced forward operator $\pmb{G}$ run faster. This is particulary relevant during the online computation of multiple posterior means for a series of measurements.  OLR is optimal for linear Gaussian Bayesian inference. For BT, there are generalisations to time-varying \cite{Sandberg2004TimeVar} and nonlinear \cite{Benner2017BilinBT,kramer2022nonlinear,kramer2022nonlinear2} systems that should be further explored in the context of data assimilation.

\subsection{TLBT for Unstable Systems}
Time-limited BT for unstable systems is applied to advection and diffusion for an instantaneous point release. For simplicity, we assume the system is one-dimensional and consider isotropic and homogeneous diffusion. For $K(z,t)$, the location- and time-dependent concentration, the process is described by
\begin{align*}
    \frac{\partial K}{\partial t} = D\frac{\partial^2K}{\partial z^2} - u\frac{\partial K}{\partial z}.
\end{align*}
The diffusion coefficient is $D =0.02$ and the mean flow is $u=0.01$. The system is discretised using finite differences for state dimensions $d = 200$ and $d = 1200$, respectively. The system matrix $\pmb{A}$ has both positive and negative eigenvalues and is unstable. The observation operator $\pmb{C}$ with $\pmb{y} = \pmb{Cx} = \begin{bmatrix}\frac{1}{d} \ldots \frac{1}{d}\end{bmatrix}\pmb{x}$ gives the mean concentration  and $d_\mathrm{out}=1$. The observation error covariance is set to $\pmb{\Gamma}_{\epsilon}= \sigma_{\mathrm{obs}}^2$, $\sigma_{\mathrm{obs}} = 0.008$. The prior covariance is chosen to be $\pmb{\Gamma}_{\mathrm{pr}} = \pmb{I}_d$ and is used as time-limited reachability Gramian, \mbox{$\pmb{P}_{\mathcal{T}}^{\text{LG}} = \pmb{\Gamma}_{\mathrm{pr}}$}. This Gramian is not prior-compatible, i.e., $\pmb{A}\pmb{\Gamma}_{\mathrm{pr}} + \pmb{\Gamma}_{\mathrm{pr}}\pmb{A}^{\mathrm{T}} = \pmb{A}+ \pmb{A}^{\mathrm{T}}$ is not negative-semidefinite, but as discussed in Sect.\@ \ref{BT-DA-pr} it is not necessary to assume prior-compatibility.

After computing the observability Gramian $\pmb{Q}_{\mathcal{T}}^{\text{LG}}$ (see Table \ref{TableGramians}), the algorithm is the same as for stable systems. To account for the instability of the system, which makes it unsuitable for long-term predictions, we use smaller end times than in the stable examples. The BT approach with infinite Gramians is not applicable to unstable sys\-tems. We, therefore, use a different approach by Qian et al.\@ \cite{Qian2021Balancing} for comparison: the Fisher matrix $\pmb{H}$ is set as the observability Gramian. No Lyapunov equation solution is required, but the computation of the full Fisher matrix is expensive. We refer to this approach as BT-H.
\begin{figure}[h]
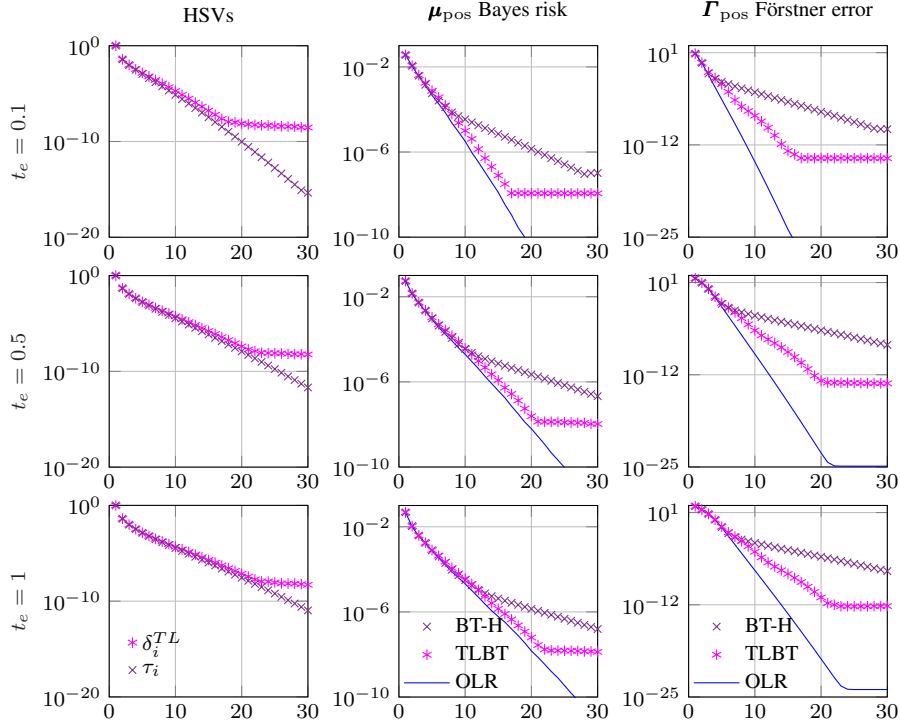

\definecolor{mycolor1}{rgb}{1.00000,0.00000,1.00000}%
\definecolor{mycolor2}{rgb}{0.49412,0.18431,0.55686}%
%
    \vspace{-3ex}
    \caption{Comparison of time-limited balanced truncation (TLBT) and balanced truncation with Fisher matrix $\pmb{H}$ as observability Gramian (BT-H) \cite{Qian2021Balancing} and the optimal low-rank approach (OLR) \cite{Spantini2015Optimal} of LG Bayesian inference for the unstable discretised advection-diffusion equation; $d = 200$. Measurements are spaced $h=0.001$ apart inside $\mathcal{T} = [0,t_e]$ for three different end times $t_e = 0.1,0.5,1$}\label{3x3 AdDiff}
\end{figure}
\begin{figure}[h]
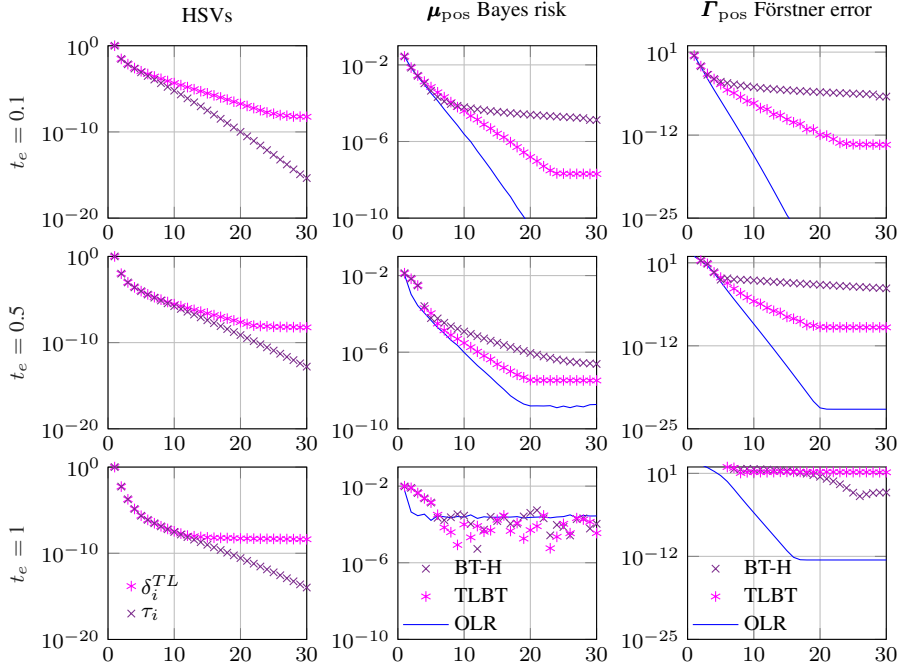

\definecolor{mycolor1}{rgb}{1.00000,0.00000,1.00000}%
\definecolor{mycolor2}{rgb}{0.49412,0.18431,0.55686}%
%    
    \vspace{-3ex}
    \caption{Model and plotted quantities as in Figure \ref{3x3 AdDiff}; here for $d = 1200$. Measurements are equispaced with \mbox{$h=0.001$} inside $\mathcal{T} = [0,t_e]$ for three different end times $t_e = 0.1,0.5,1$}\label{3x3 AdDiff_1200}
 \end{figure}

Fig.\@ \ref{3x3 AdDiff} illustrates that in our case with measurements spread quite far apart \linebreak\mbox{($h = 0.001$),} TLBT surpasses BT-H. With higher measurement frequency, the Fisher matrix approaches the time-limited observability Gramian and the results become similar. OLR ignores system stability and can be applied regardless of the system matrix $\pmb{A}$. As demonstrated in Fig.\@ \ref{3x3 AdDiff}, TLBT recovers the optimal posterior mean Bayesian risk well, but OLR outperforms TLBT for the posterior covariance approx\-imation in the Förstner metric. The error rates of TLBT (better than~$10^{-12}$) are still impressive and sufficient for applications.

For Fig.\@ \ref{3x3 AdDiff_1200} we increase the system dimension to $d=1200$. This introduces pos\-itive eigenvalues $\lambda_+$ into the unstable system matrix $\pmb{A}$, which are about an order of magnitude higher than before. The increasing instability affects the numerical proce\-dure: all posterior predictions become less reliable, especially for the largest end time $t_e=1$. This is caused by the blow-up of the values of $e^{\lambda_+\tau}$ in $e^{\pmb{A}\tau}$ for $\pmb{Q}_{\mathcal{T}}^{\text{LG}}$, already for moderate end times $t_e$. In the mean prediction, TLBT outperforms the state-of-the-art BT-H, which is not robust to instabilities since it uses exponentials of $\pmb{A}$. OLR mean prediction also fails for unstable systems and only achieves the same level of error as TLBT. These results emphasise the importance of using a time-limited ap\-proach and short observation intervals for unstable systems.

The generalisation of BT to TLBT works well, but has a main disadvantage: the numerical algorithms require the solution of a Lyapunov equation. This is only numerically feasible for (anti)stable $\pmb{A}$. In TLBT for unstable systems, the integral for $\pmb{Q}_{\mathcal{T}}^{\text{LG}}$ has to be approximated or calculated according to its definition \cite{Kuerschner2018TLBT}. To make TLBT faster and more generally applicable to large problems, further work in this direction is needed.

In this section we have demonstrated how TLBT extends the idea of Gramian-based model reduction for Bayesian inference to unstable systems. For short obser\-vation intervals TLBT significantly outperforms BT with infinite Gramians, despite the increased numerical cost of approximating the posterior statistics. 

\section{Conclusion}\label{Summary}
Model order reduction is central to data assimilation because models and data are typically high-dimensional and expensive to process. To this end, the adaptation of model reduction techniques from systems theory to data assimilation problems is essential. This work has further strengthened the link between the two research areas.

Time-limited balanced truncation (TLBT) enabled us to significantly generalise the theory of balancing Bayesian inference to systems with an unstable system matrix $\pmb{A}$ and non-compatible priors. TLBT is, thus, particularly suited to data assimilation applications such as numerical weather prediction, where short-term forecasts are required  for unstable systems \cite{Boess20114DVarunstable,BonavitaNumWeath2021}. One of the popular techniques is incremental 4D-Var and TLBT is perfectly tailored to model reduction of unstable linear tangent models in this context. The numerical efficiency of TLBT needs further investigation, e.g., by using more advanced rational Krylov methods.

Interesting research questions in data assimilation include inhomogeneous initial conditions, errors in the dynamical forward model or time-varying systems. Simi\-larities with systems theory may be instructive. Literature on such problems in the deterministic setting is available, but not extensive; see, e.g., \mbox{\cite{Beattie2016inhom,Heink2011inhom,Lang2016TimeVar,Sandberg2004TimeVar}}.

Although our approach of using TLBT in data assimilation is more general than the existing ones, it is so far limited to linear Gaussian Bayesian inference. The next step in the exploration of model order reduction for Bayesian inference is the general\-isation to nonlinear settings. The 4D-Var method could be a step in the right direction, as it already deals with nonlinear operators. Concepts from systems theory for the re\-duction of nonlinear dynamical systems may also prove useful here. This opens up a broad field of investigation at the interface between the two research communities.

\bibliography{main.bib}   % name your BibTeX data base
\bibliographystyle{spmpsci}      % mathematics and physical sciences
\end{document}